\documentclass[10pt]{amsart}
\usepackage{graphicx}
\usepackage{amsmath}
\usepackage{amssymb}
\usepackage{amsthm}
\usepackage{esint}

\newtheorem{theorem}{Theorem}[section]
\newtheorem{lemma}[theorem]{Lemma}

\theoremstyle{Corollary}
\newtheorem{cor}[theorem]{Corollary}
\newtheorem{prop}[theorem]{Proposition}
\newtheorem{definition}[theorem]{Definition}

\newtheorem{example}[theorem]{Example}

\def\R{\mathbb{R}}

\def\N{\mathbb{N}}

\numberwithin{equation}{section}

\begin{document}

\title{Some results on the weighted Yamabe problem with or without boundary}


\author{Pak Tung Ho}
\address{Department of Mathematics, Sogang University, Seoul, 04107, Korea}

\address{Korea Institute for Advanced Study, Seoul, 02455, Korea}

\email{paktungho@yahoo.com.hk}

\author{Jinwoo Shin}
\address{Korea Institute for Advanced Study, Hoegiro 85, Seoul 02455, Korea}
\email{shinjin@kias.re.kr}

\author{Zetian Yan}
\address{109 McAllister Building, Penn State University, University Park, PA 16802, USA}
\email{zxy5156@psu.edu}


\date{\today}

\begin{abstract}
Let $(M^n,g,e^{-\phi}dV_g,e^{-\phi}dA_g,m)$ be a compact smooth metric measure space with boundary with $n\geqslant 3$. In this article, we consider several Yamabe-type problems on a compact smooth metric measure space with or without boundary: uniqueness problem on the weighted Yamabe problem with boundary, characterization of the weighted Yamabe solitons with boundary and the existence of positive minimizers in the weighted Escobar quotient.

\end{abstract}

\maketitle

\section{Introduction}

Suppose $M$ is a compact, $n$-dimensional manifold
without boundary, where $n\geqslant 3$, and $g_0$ is a Riemannian metric on $M$.
As a generalization of the Uniformization Theorem, the Yamabe problem is to find a
metric conformal to $g_0$ such that its scalar curvature is constant.
This was first introduced by Yamabe \cite{Yamabe},
and was solved by Aubin \cite{Aubin0}, Trudinger \cite{Trudinger} and Schoen \cite{Schoen}.

Now consider a compact, $n$-dimensional manifold $M$
with smooth boundary $\partial M$, where $n\geqslant 3$,
and $g_0$ is a Riemannian metric on $M$.
One can still talk about the Yamabe problem
for manifold with boundary, and there are two types.
For the first type, one would like to find a conformal metric $g$
such that its scalar curvature $R_g$ is constant in $M$ and
its mean curvature $H_g$ is zero on $\partial M$.
For the second type, one would like to find a
conformal metric $g$
such that its scalar curvature $R_g$ is zero in $M$ and
its mean curvature $H_g$ is constant on $\partial M$.
These problems have studied by many authors.
See \cite{Escobar2,Escobar1,Escobar3} for example.

In this article, we consider several Yamabe-type problems on a compact smooth metric measure space with or without boundary. To explain the results of this paper, we require some terminology.

\begin{definition}
Let $(M,\partial M, g)$ be a Riemannian manifold with boundary $\partial M$ and let us denote by $dV_g$ and $dA_g$ the volume element induced by $g$ on $M$ and $\partial M$, respectively. Set a function $\phi\in C^{\infty}(M)$ and a dimensional parameter $m\in [0,\infty)$. In the case $m=0$, we require that $\phi=0$.
A smooth metric measure space (SMMS) with boundary is
a five-tuple $(M,g,e^{-\phi}dV_g,e^{-\phi}dA_g,m)$. We frequently denote a smooth metric measure space by $(M,g, v^mdV_g,v^mdA_g,m)$ where $\phi$ and $v$ are related by
$e^{-\phi}=v^m$.
\end{definition}
The \textit{weighted scalar curvature} $R^m_{\phi}$ of a smooth metric measure space with boundary $(M,g,e^{-\phi}dV_g,e^{-\phi}dA_g,m)$ is
\begin{equation}\label{defn1}
R^m_{\phi}:=R_g+2\Delta_g \phi-\frac{m+1}{m}|\nabla\phi|_g^2,
\end{equation}
where $R_g$ and $\Delta_g$ are the scalar curvature and the Laplacian associated to the metric $g$, respectively. The \textit{weighted mean curvature} is
\begin{equation}
H^m_{\phi}=H_{g}+\frac{\partial\phi}{\partial\nu_{g}},
\end{equation}
where $H_{g}$ and $\displaystyle\frac{\partial}{\partial\nu_{g}}$
are the mean curvature and the outward normal derivative with respect to $g$, respectively.

Conformal equivalence between smooth metric measure spaces are defined as follows, see \cite{Case15} for more details.
\begin{definition}\label{condef}
Smooth metric measure spaces with boundary $(M, g, e^{-\phi}dV_g,e^{-\phi}dA_g, m)$
and $(M,g_0,e^{-\phi_0}dV_{g_0},e^{-\phi_0}dA_{g_0}, m)$  are pointwise conformally equivalent if there is a function $\sigma\in C^{\infty}(M)$ such that
\begin{equation}\label{conformal}
(M, g, e^{-\phi}dV_g,e^{-\phi}dA_g, m)=(M,e^{\frac{2\sigma}{m+n-2}}g_0, e^{\frac{(m+n)\sigma}{m+n-2}}e^{-\phi_0}dV_{g_0},e^{\frac{(m+n-1)\sigma}{m+n-2}}e^{-\phi_0}dA_{g_0}, m).
\end{equation}
In the case $m=0$, this definition of conformal equivalence is reduced to the classic case.
\end{definition}
If we denote $e^{\frac{1}{2}\sigma}$ by $w$, (\ref{conformal}) is equivalent to
\begin{equation}\label{1.5}
(M, g, e^{-\phi}dV_g,e^{-\phi}dA_g, m)=(M,w^{\frac{4}{m+n-2}}g_0,w^{\frac{2(m+n)}{m+n-2}}e^{-\phi_0}dV_{g_0},w^{\frac{2(m+n-1)}{m+n-2}}e^{-\phi_0}dA_{g_0}, m),
\end{equation}
which is an alternative way to formulate the conformal equivalence of smooth metric measure spaces.

\begin{definition}
Let $(M, g, e^{-\phi}dV_g, m)$ be a smooth metric measure space. The weighted Laplacian $\Delta_{\phi}: C^{\infty}(M)\to C^{\infty}(M)$ is the operator defined as
\begin{displaymath}
\Delta_{\phi}\psi=\Delta\psi-\langle\nabla \phi, \nabla \psi\rangle_g ~~\mbox{ for any }\psi\in C^\infty(M),
\end{displaymath}
It is formally self-adjoint with respect to the measure $e^{-\phi}dV_{g}$. For more about smooth metric measure spaces, we refer the readers to \cite{Case15,Case12,Case13,Han}.
\end{definition}

\begin{definition}
Given a smooth metric measure spaces with boundary
\\ $(M, g, e^{-\phi}dV_g,e^{-\phi}dA_g, m)$, we define
\begin{equation}\label{1.7}
\begin{split}
L_{\phi}^m&=-\frac{4(n+m-1)}{n+m-2}\Delta_{\phi}+R^m_{\phi} ~~\mbox{ in }M,\\
B_{\phi}^m&=\frac{2(n+m-1)}{n+m-2}\frac{\partial}{\partial\nu_{g}}+H^m_{\phi} ~~\mbox{ on }\partial M,
\end{split}
\end{equation}
where $\nu_g$ is the outward unit normal with respect to $g$.
\end{definition}

Under the conformal change as in (\ref{1.5}),
we have the transformation laws of the weighted scalar curvature and the weighted mean curvature \cite[Proposition 1]{Posso18}:
\begin{equation}\label{1.6}
\begin{split}
R^m_{\phi}&=w^{-\frac{m+n+2}{m+n-2}}L^m_{\phi_0}w ~~\mbox{ in }M,\\
H^m_{\phi}&=w^{-\frac{n+m}{n+m-2}}B_{\phi_0}^mw~~\mbox{ on }\partial M.
\end{split}
\end{equation}
Note that our $L_{\phi}^m$ and $B_{\phi}^m$ are different than that of \cite{Posso18}
by a multiple constant. As a result, our transformation law takes the form of (\ref{1.6}).

\subsection{Uniqueness problem}

The \textit{weighted Yamabe problem} on a closed smooth metric measure space $(M,g_0,e^{-\phi_0}dV_{g_0}, m)$ is to find
another smooth metric measure space $(M,g,e^{-\phi}dV_{g}, m)$
conformal to $(M,g_0,e^{-\phi_0}dV_{g_0}, m)$ such that
its weighted scalar curvature $R^m_{{\phi}}$ is constant.
The weighted Yamabe problem in this article is different from that introduced by Case in \cite{Case15}.
See also \cite{Case12,Case19,Souza20,Posso21} for more results related to Case's weighted Yamabe problem.

Similarly, we can formulate two types of
the weighted Yamabe problem with boundary.
For the first type, one would like to find
 $(M,g,e^{-\phi}dV_{g},e^{-\phi}dA_{g}, m)$ conformal to $(M,g_0,e^{-\phi_0}dV_{g_0},e^{-\phi_0}dA_{g_0}, m)$ such that $R_{\phi}^m$ is constant in $M$
and $H_{\phi}^m$ is zero on $\partial M$.
For the second type, one would like to find
  $(M,g,e^{-\phi}dV_{g},e^{-\phi}dA_{g}, m)$ conformal to $(M,g_0,e^{-\phi_0}dV_{g_0},e^{-\phi_0}dA_{g_0}, m)$ such that $R_{\phi}^m$ is zero in $M$
and $H_{\phi}^m$ is constant on $\partial M$.

In \cite{Escobar3}, Escobar considered the uniqueness and non-uniqueness of the classic Yamabe problem with boundary. Following his argument, we study the uniqueness and non-uniqueness of the weighted Yamabe problem with boundary. Thus the precise question we will address is the following: given $(M,g,e^{-\phi}dV_{g},e^{-\phi}dA_{g}, m)$ conformal to $(M,g_0,e^{-\phi_0}dV_{g_0},e^{-\phi_0}dA_{g_0}, m)$ such that $R_{\phi}^m=R_{\phi_0}^m$ in $M$, and $H_{\phi}^m=H_{\phi_0}^m$ on $\partial M$, when we have $$(M,g,e^{-\phi}dV_{g},e^{-\phi}dA_{g}, m)=(M,g_0,e^{-\phi_0}dV_{g_0},e^{-\phi_0}dA_{g_0}, m)\mbox{ ?}$$

It follows from (\ref{1.5}) and (\ref{1.6}) that the above geometric question is equivalent to the following uniqueness question in PDEs: Suppose that $w$ is a solution of
\begin{equation}\label{2.1}
\begin{split}
-\frac{4(n+m-1)}{n+m-2}\Delta_{\phi_0}w+R^m_{\phi_0}w&=R^m_{\phi_0}w^{\frac{m+n+2}{m+n-2}}~~\mbox{ in }M,\\
\frac{2(n+m-1)}{n+m-2}\frac{\partial w}{\partial\nu_{g_0}}+H^m_{\phi_0}w&=H^m_{\phi_0}w^{\frac{n+m}{n+m-2}}~~\mbox{ on }\partial M.
\end{split}
\end{equation}
Is the function $w$ the constant function $1$ ?

First, we observe that if in the sense of \eqref{1.5}, $(M,g,e^{-\phi}dV_{g},e^{-\phi}dA_{g}, m)$ is conformal to
$(M,g_0,e^{-\phi_0}dV_{g_0},e^{-\phi_0}dA_{g_0}, m)$
such that $R^m_{\phi}=R^m_{\phi_0}=0$ in $M$ and
$H^m_{\phi}=H^m_{\phi_0}=0$ on $\partial M$,
then it follows from (\ref{1.7}) and (\ref{1.6}) that
$$\Delta_{\phi_0}w=0\mbox{ in }M~~\mbox{ and }~~\frac{\partial w}{\partial\nu_{g_0}}=0\mbox{ on }\partial M.$$
This implies that $w$ is constant, i.e.
$g=\lambda g_0$ for some positive constant $\lambda$.
Thus, from now on, we assume that  $R^m_{\phi}$ and $H^m_{\phi}$ do not vanish
simultaneously.

In the same spirit as \cite{Escobar3}, on $(M,g,e^{-\phi}dV_{g},e^{-\phi}dA_{g}, m)$, we introduce the following operators:
\begin{equation}\label{2.2}
\begin{split}
\overline{L_{\phi}^m}&=-\Delta_{\phi}-\frac{1}{n+m-1}R^m_{\phi} ~~\mbox{ in }M,\\
\overline{B_{\phi}^m}&=\frac{\partial}{\partial\nu_{g}}-\frac{1}{n+m-1}H^m_{\phi} ~~\mbox{ on }\partial M.
\end{split}
\end{equation}

Note that $(\overline{L_{\phi}^m}, \overline{B_{\phi}^m})$ is the linearization of
the problem (\ref{2.1}) at the solution $w=1$ up to resacling.

We denote $\lambda_1({L_{\phi}^m},{B_{\phi}^m})$ and $\lambda_1(\overline{L_{\phi}^m},\overline{B_{\phi}^m})$ the first eigenvalue of the
operators $({L_{\phi}^m}, {B_{\phi}^m})$ and $(\overline{L_{\phi}^m}, \overline{B_{\phi}^m})$, respectively.

Our uniqueness theorem is

\begin{theorem}\label{thm2.2}
Suppose that in the sense of \eqref{1.5}, $(M,g,e^{-\phi}dV_{g},e^{-\phi}dA_{g}, m)$ is conformal to
$(M,g_0,e^{-\phi_0}dV_{g_0},e^{-\phi_0}dA_{g_0}, m)$
with $R^m_{\phi}=R^m_{\phi_0}$ in $M$ and
$H^m_{\phi}=H^m_{\phi_0}\leq 0$ on $\partial M$. If both $\lambda_1(\overline{L_{\phi}^m},\overline{B_{\phi}^m})$
and $\lambda_1(\overline{L_{\phi_0}^m},\overline{B_{\phi_0}^m})$
are positive
or one of them is zero, then
$$(M,g,e^{-\phi}dV_{g},e^{-\phi}dA_{g}, m)=
(M,g_0,e^{-\phi_0}dV_{g_0},e^{-\phi_0}dA_{g_0}, m).$$
\end{theorem}

As a direct consequence of Theorem \ref{thm2.2}, we have the following corollary.
\begin{cor}
Suppose that in the sense of \eqref{1.5}, $(M,g,e^{-\phi}dV_{g},e^{-\phi}dA_{g}, m)$ is conformal to
$(M,g_0,e^{-\phi_0}dV_{g_0},e^{-\phi_0}dA_{g_0}, m)$
with $R^m_{\phi}=R^m_{\phi_0}\leqslant 0$ in $M$ and
$H^m_{\phi}=H^m_{\phi_0}\leqslant 0$ on $\partial M$.
Then we have
$$(M,g,e^{-\phi}dV_{g},e^{-\phi}dA_{g}, m)=
(M,g_0,e^{-\phi_0}dV_{g_0},e^{-\phi_0}dA_{g_0}, m).$$
\end{cor}

In the case that there is no uniqueness, one may consider the question of existence of metrics with the same weighted scalar and mean curvature but with one smooth metric measure space with boundary dominating the other. More precisely, we say that $[(g_1, \phi_1)]<[(g_2,\phi_2)]$ if $g_1=\varphi g_2$ and $e^{-\phi_1}=\varphi^{\frac{m}{2}}e^{-\phi_2}$, where $\varphi$ is a smooth function satisfying $0<\varphi<1$. It's natural to consider the following question: does there exist $[(g,\phi)]<[(g_0,\phi_0)]$ with the same weighted scalar and mean curvature as the one for $[(g_0,\phi_0)]$. We answer this question with the following theorem.

\begin{theorem}\label{thm3}
Suppose that $(M,g_0,e^{-\phi_0}dV_{g_0},e^{-\phi_0}dA_{g_0}, m)$  is a SMMS with boundary
such that $H_{\phi_0}^m\leq 0$.
There exists  $(M,g,e^{-\phi}dV_{g},e^{-\phi}dA_{g}, m)$ conformal to
$(M,g_0,e^{-\phi_0}dV_{g_0},e^{-\phi_0}dA_{g_0}, m)$ in the sense of \eqref{1.5}
such that $[(g,\phi)]<[(g_0,\phi_0)]$ with $R^m_{\phi}=R^m_{\phi_0}$ in $M$ and
$H^m_{\phi}=H^m_{\phi_0}\leq 0$ on $\partial M$
if and only if
$\lambda_1(L_{\phi_0}^m, B_{\phi_0}^m)<0$ and $\lambda_1(\overline{L_{\phi_0}^m}, \overline{B_{\phi_0}^m})<0$.
Furthermore, there exists at most one such metric $g$.
\end{theorem}
Moreover, we give geometric conditions that guarantee when the first eigenvalues in Theorem \ref{thm3} are negative in Section 3. We show that if $$\displaystyle\int_MR_{\phi_0}^m e^{-\phi_0}dV_{g_0}
+\int_{\partial M}H_{\phi_0}^me^{-\phi_0}dA_{g_0}\geq 0$$ then $\lambda_1(\overline{L_{\phi_0}^m}, \overline{B_{\phi_0}^m})<0$ in Proposition \ref{prop2.11} and if $$\displaystyle\int_MR_{\phi_0}^m e^{-\phi_0}dV_{g_0}
+2\int_{\partial M}H_{\phi_0}^me^{-\phi_0}dA_{g_0}\leq 0$$ then $\lambda_1(L_{\phi_0}^m, B_{\phi_0}^m)<0$ in Proposition \ref{prop2.12}.
\\

\subsection{The weighted Yamabe solitons with boundary} In this subsection, we study the weighted Yamabe solitons with boundary from two perspectives and characterize them in Theorem \ref{thm7.1} and theorem \ref{thm7.3}.

The \textit{unnormalized weighted Yamabe flow with boundary} is the evolution equation defined on $(M,g(t),e^{-\phi(t)}dV_{g(t)},e^{-\phi(t)}dA_{g(t)},m)$ given by
\begin{equation}\label{unnormalized}
  \left\{
    \begin{array}{ll}
      \displaystyle\frac{\partial}{\partial t}g(t)=-R_{\phi(t)}^mg(t) & \hbox{} \\
      \displaystyle\frac{\partial}{\partial t}\phi(t)=\frac{m}{2}R_{\phi(t)}^m & \hbox{}
    \end{array} \textrm{ in }M \textrm{ and } H_{\phi(t)}^m=0 \textrm{ on }\partial M
  \right.
\end{equation}
with $\left. g(t)\right|_{t=0}=g_0$ and $\left. \phi(t)\right|_{t=0}=\phi_0$.

From the geometric flow point of view, the weighted Yamabe solitons are the self-similar solutions to the unnormalized weighted Yamabe flow with boundary. Precisely, we say that $(M,g(t),e^{-\phi(t)}dV_{g(t)},e^{-\phi(t)}dA_{g(t)}, m)$ is a \textit{weighted Yamabe soliton with boundary} if
\begin{equation}\label{eq2.8}
       g(t)=\sigma(t)\psi_t^*(g_0)\ \ \textrm{ and }\ \ \phi(t)=\psi_t^*(\phi_0)-\frac{m}{2}\log \sigma(t)
\end{equation}
satisfies \eqref{unnormalized}, where $\sigma:[0,T)\to (0,\infty)$ is a differentiable function such that $\sigma(0)=1$, and $\psi_t:\overline{M}= M\cup \partial M\to \overline{M}$ is a 1-parameter family of diffeomorphisms in $\overline{M}$ such that $\psi_0=id_{\overline{M}}$.

The weighted Yamabe soliton can also be derived from the equation point of view. We say that a $(M,g_0,X,\phi_0,\lambda)$ is a \textit{weighted Yamabe soliton with boundary} if
\begin{equation}\label{2.12def}
  \left\{
    \begin{array}{ll}
      (\lambda-R_{\phi_0}^m)g_0=\mathcal{L}_Xg_0 &  \\
      \mathcal{L}_X\phi_0=\frac{m}{2}(R_{\phi_0}^m-\lambda) &     \end{array}
  \right.
\mbox{ in }M~~\mbox{ and }~~ H_{\phi_0}^m=0,~~X\perp \nu_{g_0}\mbox{ on }\partial M,
\end{equation}
see Section 4 for detailed discussion. 

If the vector field $X$ is a gradient vector field, i.e. $X=\frac{1}{2}\nabla f$ for some smooth function $f$ in $M$, the weighted Yamabe soliton with boudnary is called \textit{gradient} and $f$ is called \textit{the potential function}. In this case, \eqref{2.12def} can be written as
\begin{equation}\label{eq2.21}
\left\{
  \begin{array}{ll}
    \nabla^2f =(\lambda-R_{\phi_0}^m)g_0 &  \\
    \langle \nabla f,\nabla \phi_0\rangle =m(R_{\phi_0}^m-\lambda) & 
  \end{array}
\right.\textrm{ in }M \textrm{ and } H_{\phi_0}^m=0, \frac{\partial f}{\partial \nu_{g_0}}=0\textrm{ on }\partial M
\end{equation}
where $\nabla^2f$ is the Hessian of $f$ with respect to $g$.

The equivalence of (\ref{eq2.8}) and (\ref{2.12def}) is showed in Lemma \ref{sollem2.1}. Moreover, we give the characterization in the following theorems. 
\\{\bf Theorem 4.1.} Any compact weighted Yamabe soliton with boundary \eqref{eq2.8} must have constant weighted scalar curvature in $M$ and vanishing weighted mean curvature on $\partial M$.
\\{\bf Theorem 4.3.} Any compact gradient weighted Yamabe soliton with boundary \eqref{eq2.21} with $\frac{\partial \phi_0}{\partial \nu_{g_0}}=0$ must have constant weighted scalar curvature in $M$.
\\

\subsection{Existence of minimizers in the weighted Escobar quotient}
Finally, we consider the Escobar--Riemann mapping type problem on smooth metric measure spaces introduced by Posso in \cite{Posso18} and prove the conjecture left there.

The Escobar--Riemann mapping problem on a manifold with boundary $(M^n,\partial M,g)$ is concerned with finding a metric $g$ with vanishing scalar curvature in $M$ and constant mean curvature on $\partial M$, within the conformal class of $g$. This problem in the Euclidean half-space reduces to finding the minimizers in the sharp Trace Gagliardo--Nirenberg--Sobolev inequality. To present the Trace Gagliardo--Nirenberg--Sobolev inequality, we denote the half-space and its boundary by
\begin{displaymath}
    \R^n_{+}=\left\{(x,t):x\in \R^{n-1}, t\geqslant 0\right\} \mbox{~~and~~} \partial\R^n_+=\left\{(x,0)\in \R^n :x\in \R^{n-1}\right\}.
\end{displaymath}
We identify $\partial\R^n_+$ with $\R^{n-1}$ whenever necessary.
\begin{theorem}[\cite{BCFG}]
Fix $m\geqslant 0$. For all $w\in W^{1,2}(\R^n_+)\cap L^{\frac{2(m+n-1)}{m+n-2}}(\R^n_+)$ it holds that
\begin{equation}\label{Trace GNS}
        \Lambda_{m,n}\left(\int_{\partial\R^n_+}|w|^{\frac{2(m+n-1)}{m+n-2}}\right)^{\frac{2m+n-2}{m+n-1}}\leqslant \left(\int_{\R^n_+}|\nabla w|^2\right)\left(\int_{\R^n_+} |w|^{\frac{2(m+n-1)}{m+n-2}}\right)^{\frac{m}{m+n-1}}
\end{equation}
where the constant $\Lambda_{m,n}$ is given by
\begin{equation}
    \Lambda_{m,n}=(m+n-2)^2\left(\frac{{\rm Vol}(\mathbb{S}^{2m+n-1})^{\frac{1}{2m+n-1}}}{2(2m+n-2)}\right)^{\frac{2m+n-1}{m+n-1}}\left(\frac{\Gamma(2m+n-1)}{\pi^m\Gamma(m+n-1)}\right)^{\frac{1}{m+n-1}}
\end{equation}
and ${\rm Vol}(\mathbb{S}^{2m+n-1})$ is the volume of the $2m+n-1$ dimensional unit sphere. Moreover, equality holds if and only if $w$ is a constant multiple of the function $w_{\epsilon,x_0}$ defined on $\R^n_+$ by
\begin{equation}
    w_{\epsilon,x_0}(x,t):=\left(\frac{2\epsilon}{(\epsilon+t)^2+|x-x_0|^2}\right)^{\frac{m+n-2}{2}}
\end{equation}
where $\epsilon>0$ and $x_0\in \R^{n-1}$.
\end{theorem}

In \cite{Posso18}, the author introduced an analogue of (\ref{Trace GNS}) on smooth metric measure spaces with boundary called the weighted Escobar quotient.

\begin{definition}
Given a smooth metric measure spaces with boundary
$(M, g, e^{-\phi}dV_g,e^{-\phi}dA_g, m)$, the weighted Escobar quotient ${\mathcal{Q}}: W^{1,2}(M,e^{-\phi}dV_{g})\to \R$ is defined by
\begin{equation}\label{quotient}
{\mathcal{Q}}(w)=\frac{\int_MwL_{\phi}^m(w)e^{-\phi}dV_{g}+\int_M wB_{\phi}^m(w)e^{-\phi}dA_{g}}{\left(\int_{\partial M}|w|^{\frac{2(n+m-1)}{n+m-2}}e^{-\phi}dA_{g}\right)^{\frac{n+2m-2}{n+m-1}}\left(\int_M |w|^{\frac{2(n+m-1)}{n+m-2}}e^{-\frac{m-1}{m}\phi}dV_g\right)^{-\frac{m}{n+m-1}}}.
\end{equation}
\end{definition}
The {\textit{weighted Escobar constant}} $\Lambda[M, g, e^{-\phi}dV_g,e^{-\phi}dA_g, m]\in\R\cup \{-\infty\}$ is defined by
\begin{equation}
    \Lambda:=\Lambda[M, g, e^{-\phi}dV_g,e^{-\phi}dA_g, m]=\inf \left\{{\mathcal{Q}}(w): w\in W^{1,2}(M,e^{-\phi}dV_{g})\right\}.
\end{equation}
If $m=0$, the quotient (\ref{quotient}) conincides with the Sobolev quotient considered by Escobar in the Escobar-Riemann mapping problem.

In \cite{Posso18}, Posso proved that the infimum can be achieved by a positive function when $\Lambda<0$ and $\Lambda[M, g, e^{-\phi}dV_g,e^{-\phi}dA_g, m] \leq \Lambda_{m,n}$ on general smooth metric measure spaces with boundary. Moreover, he conjectured the existence of positive minimizers of ${\mathcal{Q}}(w)$ when $\Lambda<\Lambda_{m,n}$. 

In the rest of this article, we prove the existencce of a positive minimizer of the weighted Escobar constant in the subcritical case; i.e. $\Lambda<\Lambda_{m,n}$. 
\begin{theorem}\label{main}
Let $(M, g, e^{-\phi}dV_g,e^{-\phi}dA_g, m)$ be a smooth metric measure spaces with boundary, $m\geqslant 0$ and $\Lambda<\Lambda_{m,n}$. Then there exists a positive function $w\in C^{\infty}(M)$ such that
\begin{displaymath}
    {\mathcal{Q}}(w)=\Lambda[M, g, e^{-\phi}dV_g,e^{-\phi}dA_g, m].
\end{displaymath}
\end{theorem}

The key ideas in the proof of Theorem \ref{main} are the following. First, we adapt the argument in \cite{DHR} to complete the blow-up analysis on the weighted Escobar quotient; see Theorem \ref{blowup} for a precise statement. Second, in the same spirit as Theorem 1.1 in \cite{Yan}, we prove Aubin-type almost trace GNS
inequalities on general smooth metric measure spaces with $\Lambda \geq 0$. The proof of Theorem \ref{main} follows from these results immediately.

This paper is organized as follows.

In Section 2, in the same spirit as \cite{Escobar3}, we establish the uniqueness result-Theorem \ref{thm2.2}.

In Section 3 we prove Theorem \ref{thm3}. The proof consists of two parts. In the first part we establish the existence of a smaller smooth metric measure space with the same weighted scalar and mean curvature by the method of upper and lower solutions. In the second part we show the uniqueness of the solution we found. In order to do that we use degree theory, and compute the Leray--Schauder degree of the possible solutions. Moreover, we give two geometric criteria in Proposition \ref{prop2.11} and Proposition \ref{prop2.12}.

In Section 4 we prove the equivalence between (\ref{eq2.8}) and (\ref{2.12def}) and characterize the weighted Yamabe solitons with boundary defined in (\ref{eq2.8}) and (\ref{eq2.21}).

In Section 5 we mimic the proof in \cite{DHR} to characterize Palais--Smale sequences of the weighted Escobar quotient.

In Section 6 we establish Aubin-type almost trace GNS inequalities and complete the proof of Theorem \ref{main}.

\section{A uniqueness theorem}
The purpose of this section is to show the validity of Theorem \ref{thm2.2}. In order to do that we need the following lemmas.

\begin{lemma}
Suppose that in the sense of \eqref{1.5}, $(M,g,e^{-\phi}dV_{g},e^{-\phi}dA_{g}, m)$ is conformal to
$(M,g_0,e^{-\phi_0}dV_{g_0},e^{-\phi_0}dA_{g_0}, m)$
with $R^m_{\phi}=R^m_{\phi_0}<0$ in $M$ and
$H^m_{\phi}=H^m_{\phi_0}<0$ on $\partial M$.
Then we have
$$(M,g,e^{-\phi}dV_{g},e^{-\phi}dA_{g}, m)=
(M,g_0,e^{-\phi_0}dV_{g_0},e^{-\phi_0}dA_{g_0}, m).$$
\end{lemma}
\begin{proof}

Consider $x_0\in\overline{M}$ such that $ w(x_0)=\max_{\overline{M}}w$.
There are two cases to be considered.\\
\textit{Case 1.} If $x_0\in M$, then $\Delta_{\phi_0} w\leq 0$ at $x_0$.
This together with the first equation of (\ref{2.1}) implies that
$$R^m_{\phi_0}\big(w^{\frac{m+n+2}{m+n-2}}-w\big)\geq 0~~\mbox{ at }x_0.$$
Since $R^m_{\phi_0}<0$ by assumption, we have
$w^{\frac{m+n+2}{m+n-2}}-w\leq 0$ at $x_0$,
which gives
$w(x_0)=\max_{\overline{M}}w\leq 1$.\\
\textit{Case 2.} If $x_0\in \partial M$, then $\displaystyle\frac{\partial w}{\partial\nu_{g_0}}\geq 0$ at $x_0$.
This together with  the second equation of (\ref{2.1}) implies that
$$H^m_{\phi_0}\big(w^{\frac{m+n}{m+n-2}}-w\big)\geq 0~~\mbox{ at }x_0.$$
By assumption, $H^m_{\phi_0}<0$. Hence, we have
$w^{\frac{m+n}{m+n-2}}-w\leq 0$ at $x_0$,
which gives
$w(x_0)=\max_{\overline{M}}w\leq 1$.

In both case, we have $\max_{\overline{M}}w\leq 1$.
Similarly, by considering $x_1\in\overline{M}$ with $w(x_1)=\min_{\overline{M}}w$,
we can show that $\min_{\overline{M}}w\geq 1$.
Thus we have $w\equiv 1$, which proves the assertion.
\end{proof}

\begin{lemma}\label{lem2.3}
Suppose that in the sense of \eqref{1.5}, $(M,g,e^{-\phi}dV_{g},e^{-\phi}dA_{g}, m)$ is conformal to
$(M,g_0,e^{-\phi_0}dV_{g_0},e^{-\phi_0}dA_{g_0}, m)$
with $R^m_{\phi}=R^m_{\phi_0}$ in $M$ and
$H^m_{\phi}=H^m_{\phi_0}\leq 0$ on $\partial M$.\\
(i) If $\lambda_1(\overline{L_{\phi_0}^m},\overline{B_{\phi_0}^m})>0$, then either $g=g_0$ or $[(g,\phi)]>[(g_0,\phi_0)]$.\\
(ii) If $\lambda_1(\overline{L_{\phi_0}^m},\overline{B_{\phi_0}^m})=0$, then $g=g_0$.
\end{lemma}
\begin{proof}
It follows from the assumption that (\ref{2.1}) holds.
Let $f=w^{-\frac{4}{n+m-2}}-1$.
By (\ref{2.1}), a direct computation yields
\begin{equation}\label{2.3}
\begin{split}
\Delta_{\phi_0}f+\frac{1}{n+m-1}R_{\phi_0}^mf&=\frac{4(n+m+2)}{(n+m-2)^2}w^{-\frac{2(n+m)}{n+m-2}}|\nabla_{g_0}w|^2~~\mbox{ in }M,\\
\frac{\partial f}{\partial\nu_{g_0}}&=\frac{2}{n+m-1}\frac{H_{\phi_0}^mf}{(w^{\frac{2}{n+m-2}}+1)}~~\mbox{ on }\partial M.
\end{split}
\end{equation}
Let $\psi>0$ be the eigenfunction to $\lambda_1(\overline{L_{\phi_0}^m},\overline{B_{\phi_0}^m})$, i.e.
\begin{equation}\label{2.4}
\begin{split}
-\Delta_{\phi_0}\psi-\frac{1}{n+m-1}R^m_{\phi_0}\psi&=\lambda_1(\overline{L_{\phi_0}^m},\overline{B_{\phi_0}^m})\psi ~~\mbox{ in }M,\\
\frac{\partial\psi}{\partial\nu_{g_0}}-\frac{1}{n+m-1}H^m_{\phi_0}\psi&=0 ~~\mbox{ on }\partial M.
\end{split}
\end{equation}
Set $v=f/\psi$. Then $f=v\psi$ and
$$\Delta_{\phi_0}f=\psi\Delta_{\phi_0} v+v\Delta_{\phi_0}\psi+2\langle\nabla_{g_0}v,\nabla_{g_0}\psi\rangle.$$
This together with (\ref{2.3}) and (\ref{2.4}) gives
\begin{equation}\label{2.5}
\begin{split}
\Delta_{\phi_0}v &=-\frac{v}{\psi}\Delta_{\phi_0}\psi-\frac{2}{\psi}\langle\nabla_{g_0}v,\nabla_{g_0}\psi\rangle+\frac{1}{\psi}\Delta_{\phi_0}f\\
&=-\frac{v}{\psi}\left(-\frac{1}{n+m-1}R^m_{\phi_0}\psi-\lambda_1(\overline{L_{\phi_0}^m},\overline{B_{\phi_0}^m})\psi\right)
-\frac{2}{\psi}\langle\nabla_{g_0}v,\nabla_{g_0}\psi\rangle\\
&\hspace{4mm}+\frac{1}{\psi}\left(-\frac{1}{n+m-1}R_{\phi_0}^mf+\frac{4(n+m+2)}{(n+m-2)^2}w^{-\frac{2(n+m)}{n+m-2}}|\nabla_{g_0}w|^2\right)\\
&= \lambda_1(\overline{L_{\phi_0}^m},\overline{B_{\phi_0}^m})v-\frac{2}{\psi}\langle\nabla_{g_0}v,\nabla_{g_0}\psi\rangle\\
&\hspace{4mm}
+\frac{4(n+m+2)}{\psi(n+m-2)^2}w^{-\frac{2(n+m)}{n+m-2}}|\nabla_{g_0}w|^2\\
&\geq  \lambda_1(\overline{L_{\phi_0}^m},\overline{B_{\phi_0}^m})v-\frac{2}{\psi}\langle\nabla_{g_0}v,\nabla_{g_0}\psi\rangle~~\mbox{ in }M,
\end{split}
\end{equation}
and
\begin{equation}\label{2.6}
\frac{\partial v}{\partial\nu_{g_0}}
=\frac{1}{\psi}\frac{\partial f}{\partial\nu_{g_0}}-\frac{f}{\psi^2}\frac{\partial \psi}{\partial\nu_{g_0}} =\frac{H_{\phi_0}^mv}{n+m-1}\left(\frac{1-w^{\frac{2}{n+m-2}}}{1+w^{\frac{2}{n+m-2}}}\right)~~\mbox{ on }\partial M.
\end{equation}
Let $v(x_0)=\max_{\overline{M}}v$. There are two cases to be considered.\\
\textit{Case 1.} If $x_0\in M$, then $\Delta_{\phi_0} v\leq 0$ and $\nabla_{g_0} v=0$ at $x_0$.
By (\ref{2.5}) and the assumption that $\lambda_1(\overline{L_{\phi_0}^m},\overline{B_{\phi_0}^m})>0$,
we have $v(x_0)\leq 0$.

If $v(x_0)=0$, then the strong maximum principle implies that $v\equiv 0$, which
implies that $f\equiv 0$, or equivalently, $w\equiv 1$.

If $v(x_0)<0$, then $v<0$, which implies that $f<0$, or equivalently, $w>1$.\\
\textit{Case 2.} If $x_0\in\partial M$, then
either $\displaystyle \frac{\partial v}{\partial\nu_{g_0}}(x_0)>0$
or $v$ is constant.
In the latter case, the constant is either zero or negative by (\ref{2.5}).
Hence, we have $w\equiv 1$ or $w>1$.
The case $\displaystyle \frac{\partial v}{\partial\nu_{g_0}}(x_0)>0$
is ruled out. To see this, note that
$\displaystyle \frac{\partial v}{\partial\nu_{g_0}}(x_0)>0$
implies that
$$\frac{H_{\phi_0}^mv}{n+m-1}\left(\frac{1-w^{\frac{2}{n+m-2}}}{1+w^{\frac{2}{n+m-2}}}\right)>0~~\mbox{ at }x_0.$$
Since $v=\displaystyle\frac{w^{-\frac{4}{n+m-2}}-1}{\psi}$, this implies that
$$\frac{H_{\phi_0}^m}{n+m-1}\left(\frac{w^{-\frac{4}{n+m-2}}-1}{\psi}\right)\left(\frac{1-w^{\frac{2}{n+m-2}}}{1+w^{\frac{2}{n+m-2}}}\right)>0~~\mbox{ at }x_0,$$
which gives
$H_{\phi_0}^m(1-w^{\frac{2}{n+m-2}})^2>0$ at $x_0$.
But this contradicts to the assumption that $H_{\phi_0}^m\leq 0$.

Now let us assume that $\lambda_1(\overline{L_{\phi_0}^m},\overline{B_{\phi_0}^m})=0$.
In this case, the maximum principle implies that $x_0\in\partial M$.
Therefore, either $\displaystyle \frac{\partial v}{\partial\nu_{g_0}}(x_0)>0$
or $v$ is constant. Since we have shown that it is impossible for $\displaystyle \frac{\partial v}{\partial\nu_{g_0}}(x_0)>0$, $v$ must be constant.
Hence, it follows from (\ref{2.5}) that
$|\nabla_{g_0}w|^2=0$ in $M$, which implies that
$w$ is constant.
Since we have assumed that
$R^m_{\phi}=R^m_{\phi_0}$ and $H^m_{\phi}=H^m_{\phi_0}$ do not vanish
simultaneously, it follows from
(\ref{2.1}) that $w\equiv 1$.
\end{proof}

\begin{proof}[Proof of Theorem \ref{thm2.2}]
If either  $\lambda_1(\overline{L_{\phi}^m},\overline{B_{\phi}^m})$
or $\lambda_1(\overline{L_{\phi_0}^m},\overline{B_{\phi_0}^m})$ is zero,
it follows from Lemma \ref{lem2.3} that $g=g_0$.

On the other hand, if $\lambda_1(\overline{L_{\phi}^m},\overline{B_{\phi}^m})$
is positive, then $g_0=g$ or $g_0>g$ by Lemma \ref{lem2.3}.
If $\lambda_1(\overline{L_{\phi_0}^m},\overline{B_{\phi_0}^m})$ is positive,
then $g=g_0$ or $g>g_0$ by Lemma \ref{lem2.3}.
Therefore, if both $\lambda_1(\overline{L_{\phi}^m},\overline{B_{\phi}^m})$
and $\lambda_1(\overline{L_{\phi_0}^m},\overline{B_{\phi_0}^m})$ are positive,
the only possibility is
that $g=g_0$. This proves Theorem \ref{thm2.2}.
\end{proof}

\section{Existence and uniqueness of smaller metrics}
We prove Theorem \ref{thm3} in this section. First, we need the following preliminary results.

\begin{lemma}\label{lem2.4}
If there exists $[(g,\phi)]<[(g_0,\phi_0)]$ such that
$R^m_{\phi}=R^m_{\phi_0}$ in $M$ and
$H^m_{\phi}=H^m_{\phi_0}\leq 0$ on $\partial M$,
then $\lambda_1(L_{\phi_0}^m,B_{\phi_0}^m)<0$.
\end{lemma}
\begin{proof}
If we write $g=w^{\frac{4}{n+m-2}}g_0$,
then $w<1$ is a solution to (\ref{2.1}).
By (\ref{2.1})
and integration by parts, we obtain
\begin{equation}\label{2.7}
\begin{split}
E(w)&=\int_M\left(\frac{4(n+m-1)}{n+m-2}|\nabla_{g_0}w|^2+R_{\phi_0}^mw^2\right)e^{-\phi_0}dV_{g_0}+2\int_{\partial M} H_{\phi_0}^mw^2e^{-\phi_0}dA_{g_0}\\
&=\int_Mw\left(-\frac{4(n+m-1)}{n+m-2}\Delta_{\phi_0} w+R_{\phi_0}^mw\right)e^{-\phi_0}dV_{g_0}\\
&\hspace{4mm}+2\int_{\partial M}w\left(\frac{2(n+m-1)}{n+m-2}\frac{\partial w}{\partial\nu_{g_0}}+ H_{\phi_0}^mw\right)e^{-\phi_0}dA_{g_0}\\
&=\int_MR_{\phi_0}^mw^{\frac{2(n+m)}{n+m-2}}e^{-\phi_0}dV_{g_0}
+2\int_{\partial M} H_{\phi_0}^mw^{\frac{2(n+m-1)}{n+m-2}}e^{-\phi_0}dA_{g_0}.
\end{split}
\end{equation}
It follows from (\ref{2.1}) and integration by parts that
\begin{equation*}
\begin{split}
&\int_MR_{\phi_0}^mw^{\frac{2(n+m)}{n+m-2}}e^{-\phi_0}dV_{g_0}\\
&=\frac{4(n+m-1)}{n+m-2}\int_M\frac{w^{\frac{n+m+2}{n+m-2}}}{1-w^{\frac{4}{n+m-2}}}\Delta_{\phi_0}we^{-\phi_0}dV_{g_0}\\
&=-\frac{4(n+m-1)}{n+m-2}\int_M\left\langle\nabla_{g_0}\left(\frac{w^{\frac{n+m+2}{n+m-2}}}{1-w^{\frac{4}{n+m-2}}}\right),\nabla_{g_0}w\right\rangle e^{-\phi_0}dV_{g_0}\\
&\hspace{4mm}
+\frac{4(n+m-1)}{n+m-2}\int_{\partial M}\frac{w^{\frac{n+m+2}{n+m-2}}}{1-w^{\frac{4}{n+m-2}}}\frac{\partial w}{\partial\nu_{g_0}}e^{-\phi_0}dA_{g_0}\\
&=-\frac{4(n+m-1)(n+m+2)}{(n+m-2)^2}\int_M\frac{w^{\frac{4}{n+m-2}}}{1-w^{\frac{4}{n+m-2}}}|\nabla_{\phi_0}w|^2 e^{-\phi_0}dV_{g_0}\\
&\hspace{4mm}-\frac{16(n+m-1)}{(n+m-2)^2}\int_M\frac{w^{\frac{8}{n+m-2}}}{(1-w^{\frac{4}{n+m-2}})^2}|\nabla_{\phi_0}w|^2 e^{-\phi_0}dV_{g_0}
-2\int_{\partial M}H_{\phi_0}^m \frac{w^{\frac{2(n+m)}{n+m-2}}}{1+w^{\frac{2}{n+m-2}}}e^{-\phi_0}dA_{g_0}.
\end{split}
\end{equation*}
Combining this with (\ref{2.7}) yields
\begin{equation*}
\begin{split}
E(w)&=-\frac{4(n+m-1)(n+m+2)}{(n+m-2)^2}\int_M\frac{w^{\frac{4}{n+m-2}}}{1-w^{\frac{4}{n+m-2}}}|\nabla_{\phi_0}w|^2 e^{-\phi_0}dV_{g_0}\\
&\hspace{4mm}-\frac{16(n+m-1)}{(n+m-2)^2}\int_M\frac{w^{\frac{8}{n+m-2}}}{(1-w^{\frac{4}{n+m-2}})^2}|\nabla_{\phi_0}w|^2 e^{-\phi_0}dV_{g_0}\\
&\hspace{4mm}
+2\int_{\partial M} H_{\phi_0}^m\frac{w^{\frac{2(n+m-1)}{n+m-2}}}{1+w^{\frac{2}{n+m-2}}}e^{-\phi_0}dA_{g_0}.
\end{split}
\end{equation*}
Since $w<1$ and $H_{\phi_0}^m\leq 0$ by assumption, we have $E(w)<0$.
From the variational characterization of $\lambda_1(L_{\phi_0}^m,B_{\phi_0}^m)$,
$$\lambda_1(L_{\phi_0}^m,B_{\phi_0}^m)=\inf_{0\not\equiv u\in C^\infty(M)}\frac{E(u)}{\int_M |u|^2e^{-\phi_0}dV_{g_0}},$$
we find that $\lambda_1(L_{\phi_0}^m,B_{\phi_0}^m)<0$.
\end{proof}

\begin{lemma}\label{lem2.5}
If there exists $[(g,\phi)]<[(g_0,\phi_0)]$ such that
$R^m_{\phi}=R^m_{\phi_0}$ in $M$ and
$H^m_{\phi}=H^m_{\phi_0}\leq0$ on $\partial M$,
then $\lambda_1(\overline{L_{\phi_0}^m},\overline{B_{\phi_0}^m})<0$.
\end{lemma}
\begin{proof}
By the assumption $[(g,\phi)]<[(g_0,\phi_0)]$, we have $w<1$ which implies
\begin{equation}\label{mistake}
  \frac{2}{1+w^\frac{2}{n+m-2}}>1.
\end{equation}

It follows from \eqref{2.3}, \eqref{mistake}, and $H_g\leq 0$ that the function $f=w^{-\frac{4}{n+m-2}}-1>0$ satisfies
\begin{equation}\label{2.8}
\begin{split}
\Delta_{\phi_0}f+\frac{1}{n+m-1}R_{\phi_0}^mf>0&~~\mbox{ in }M,\\
\frac{\partial f}{\partial\nu_{g_0}}\leq\frac{1}{n+m-1}H_{\phi_0}^mf&~~\mbox{ on }\partial M.
\end{split}
\end{equation}
Multiplying the first equation by $f$ and integrating by parts, we obtain
\begin{equation*}
\begin{split}
0&<\int_M\left(f\Delta_{\phi_0}f+\frac{1}{n+m-1}R_{\phi_0}^mf^2\right)e^{-\phi_0}dV_{g_0}\\
&=\int_M\left(-|\nabla_{\phi_0}f|^2+\frac{1}{n+m-1}R_{\phi_0}^mf^2\right)e^{-\phi_0}dV_{g_0}
+\int_{\partial M}f\frac{\partial f}{\partial\nu_{g_0}}e^{-\phi_0}dA_{g_0}\\
&\leq\int_M\left(-|\nabla_{\phi_0}f|^2+\frac{1}{n+m-1}R_{\phi_0}^mf^2\right)e^{-\phi_0}dV_{g_0}\\
&\hspace{4mm}
+\frac{1}{n+m-1}\int_{\partial M}H_{\phi_0}^m f^2e^{-\phi_0}dA_{g_0}.
\end{split}
\end{equation*}
Therefore, from the variational characterization of $\lambda_1(\overline{L_{\phi_0}^m},\overline{B_{\phi_0}^m})$,
\begin{equation*}
\begin{split}
&\lambda_1(\overline{L_{\phi_0}^m},\overline{B_{\phi_0}^m})\\
&=
\inf_{0\not\equiv u\in C^\infty(M)}\frac{\int_M\left(|\nabla_{\phi_0}f|^2-\frac{R_{\phi_0}^m}{n+m-1}f^2\right)e^{-\phi_0}dV_{g_0}
-\frac{1}{n+m-1}\int_{\partial M}H_{\phi_0}^m u^2e^{-\phi_0}dA_{g_0}}{\int_M |u|^2e^{-\phi_0}dV_{g_0}},
\end{split}
\end{equation*}
we find that $\lambda_1(\overline{L_{\phi_0}^m},\overline{B_{\phi_0}^m})<0$.
\end{proof}

\begin{lemma}\label{lem2.6}
Suppose that $(M,g_0,e^{-\phi_0}dV_{g_0},e^{-\phi_0}dA_{g_0}, m)$
is a compact smooth metric measure space with boundary
such that $H_{\phi_0}^m\leq 0$, $\lambda_1(L_{\phi_0}^m,B_{\phi_0}^m)<0$
and
$\lambda_1(\overline{L_{\phi_0}^m},\overline{B_{\phi_0}^m})<0$.
There exists $[(g,\phi)]<[(g_0,\phi_0)]$ such that
$R^m_{\phi}=R^m_{\phi_0}$ in $M$ and
$H^m_{\phi}=H^m_{\phi_0}$ on $\partial M$.
\end{lemma}
\begin{proof}
Let $\varphi_1>0$ be the first eigenfunction of $(L_{\phi_0}^m,B_{\phi_0}^m)$, i.e.
\begin{equation}\label{2.9}
\begin{split}
-\frac{4(n+m-1)}{n+m-2}\Delta_{\phi_0}\varphi_1+R^m_{\phi_0}\varphi_1&=\lambda_1(L_{\phi_0}^m,B_{\phi_0}^m)\varphi_1 ~~\mbox{ in }M,\\
\frac{2(n+m-1)}{n+m-2}\frac{\partial\varphi_1}{\partial\nu_{g_0}}+H^m_{\phi_0}\varphi_1&=0 ~~\mbox{ on }\partial M.
\end{split}
\end{equation}
By rescaling, we may assume that $\max_{\overline{M}}\varphi_1=1$.
Define $u_0=\epsilon\varphi_1^\alpha$, where $\alpha$ is close to $1$ and $\epsilon$ is small.
Then $u_0$ is a lower solution to (\ref{2.1}), because
\begin{equation*}
\begin{split}
&\frac{4(n+m-1)}{n+m-2}\Delta_{\phi_0}u_0-R^m_{\phi_0}u_0+R^m_{\phi_0}u_0^{\frac{m+n+2}{m+n-2}}\\
&=\frac{4(n+m-1)}{n+m-2}\epsilon\big(\alpha\varphi_1^{\alpha-1}\Delta_{\phi_0}\varphi_1
+\alpha(\alpha-1)\varphi_1^{\alpha-2}|\nabla_{g_0}\varphi_1|^2\big)-R^m_{\phi_0}\epsilon\varphi_1^\alpha+R^m_{\phi_0}(\epsilon\varphi_1^\alpha)^{\frac{m+n+2}{m+n-2}}\\
&=\epsilon\varphi_1^\alpha\left[\frac{4(n+m-1)}{n+m-2}\alpha(\alpha-1)\frac{|\nabla_{g_0}\varphi_1|^2}{\varphi_1^2}-\alpha\lambda_1(L_{\phi_0}^m,B_{\phi_0}^m)
+R^m_{\phi_0}(\alpha-1)+R^m_{\phi_0}(\epsilon\varphi_1^\alpha)^{\frac{4}{m+n-2}}\right]
\end{split}
\end{equation*}
in $M$, where we have used the first equation in (\ref{2.9}).  Thus, if $\alpha$ is close to $1$ and $\epsilon$ is small enough, we have
$$\frac{4(n+m-1)}{n+m-2}\Delta_{\phi_0}u_0-R^m_{\phi_0}u_0+R^m_{\phi_0}u_0^{\frac{m+n+2}{m+n-2}}>0~~\mbox{ in }M$$
since $\lambda_1(L_{\phi_0}^m,B_{\phi_0}^m)<0$.
By the second equation in (\ref{2.9}), we also have
\begin{equation*}
\begin{split}
&\frac{2(n+m-1)}{n+m-2}\frac{\partial u_0}{\partial\nu_{g_0}}+H^m_{\phi_0}u_0-H^m_{\phi_0}u_0^{\frac{n+m}{n+m-2}}\\
&=\frac{2(n+m-1)}{n+m-2}\epsilon\alpha\varphi_1^{\alpha-1}\frac{\partial \varphi_1}{\partial\nu_{g_0}}+H^m_{\phi_0}\epsilon\varphi_1^\alpha-H^m_{\phi_0}(\epsilon\varphi_1^\alpha)^{\frac{n+m}{n+m-2}}\\
&=\epsilon\varphi_1^{\alpha}H^m_{\phi_0}\left[1-\alpha-(\epsilon\varphi_1^\alpha)^{\frac{2}{n+m-2}}\right]
\end{split}
\end{equation*}
on $\partial M$. If we choose $\alpha=1-\epsilon^{\frac{2}{n+m-2}}$, we can deduce that the last
expression is nonpositive, i.e.
$$\frac{2(n+m-1)}{n+m-2}\frac{\partial u_0}{\partial\nu_{g_0}}+H^m_{\phi_0}u_0-H^m_{\phi_0}u_0^{\frac{n+m}{n+m-2}}\leq 0~~\mbox{ on }\partial M.$$
Hence, $u_0$ is a lower solution to (\ref{2.1}).

To construct an upper solution, we let $f_1>0$ be the first eigenfunction to
problem (\ref{2.4}). By rescaling, we may assume that $\max_{\overline{M}}f_1=1$.
Define $w=1-\delta f_1$, where $\delta>0$.
We compute
\begin{equation*}
\begin{split}
&\frac{4(n+m-1)}{n+m-2}\Delta_{\phi_0}w-R^m_{\phi_0}w+R^m_{\phi_0}w^{\frac{m+n+2}{m+n-2}}\\
&=-\delta\frac{4(n+m-1)}{n+m-2}\Delta_{\phi_0}f_1-R^m_{\phi_0}(1-\delta f_1)
+R^m_{\phi_0}(1-\delta f_1)^{\frac{m+n+2}{m+n-2}}\\
&=\delta f_1\Bigg\{\frac{4(n+m-1)}{n+m-2}\lambda_1(\overline{L_{\phi_0}^m},\overline{B_{\phi_0}^m})\\
&\hspace{8mm}
+R_{\phi_0}^m\Bigg[\frac{4}{n+m-2}+(1-\delta f_1)\frac{\big(-1+(1-\delta f_1)^{\frac{4}{n+m-2}}\big)}{\delta f_1}\Bigg]\Bigg\}.
\end{split}
\end{equation*}
Since
$$\lim_{\delta\to 0^+}\frac{1-(1-\delta f_1)^{\frac{4}{n+m-2}}}{\delta f_1}=\frac{4}{n+m-2},$$
we get that for $\delta$ small enough
$$\frac{4(n+m-1)}{n+m-2}\Delta_{\phi_0}w-R^m_{\phi_0}w+R^m_{\phi_0}w^{\frac{m+n+2}{m+n-2}}<0~~\mbox{ in }M.$$
Now, on $\partial M$, we have
\begin{equation*}
\begin{split}
&\frac{2(n+m-1)}{n+m-2}\frac{\partial w}{\partial\nu_{g_0}}+H^m_{\phi_0}w-H^m_{\phi_0}w^{\frac{n+m}{n+m-2}}\\
&=\delta f_1 H_{\phi_0}^m\left[-\frac{2}{n+m-2}+(1-\delta f_1)\left(\frac{1-(1-\delta f_1)^{\frac{2}{n+m-2}}}{\delta f_1}\right)\right]
\end{split}
\end{equation*}
Observe that for $x>0$ small enough
$$(1-x)\left(\frac{1-(1-x)^{\frac{2}{n+m-2}}}{x}\right)\leq \frac{2}{n+m-2}.$$
Hence, taking $\delta$ small enough, we obtain
$$\frac{2(n+m-1)}{n+m-2}\frac{\partial w}{\partial\nu_{g_0}}+H^m_{\phi_0}w-H^m_{\phi_0}w^{\frac{n+m}{n+m-2}}\geq 0~~\mbox{ on }\partial M.$$
Hence, $w$ is an upper solution to (\ref{2.1}).

The method of upper and lower solutions then implies that,
by taking $\epsilon$ and $\delta$ small enough and hence $\alpha<1$ but close to $1$,
we can find a smooth solution $u$ to (\ref{2.1}) with
$R_{\phi}^m=R_{\phi_0}^m$ in $M$ and $H_{\phi}^m=H_{\phi_0}^m$ on $\partial M$
satisfying
$$0<u_0<u<w<1.$$
As a result, we get the required metric $g=u^{\frac{4}{n+m-2}}g_0$.
\end{proof}

Consider the operator $T:C^\infty(M)\to C^\infty(M)$ defined as $T(\varphi)=\psi$
where $\psi$ is the unique solution to the boundary value problem
\begin{equation*}
\begin{split}
\Delta_{\phi_0}\psi-\gamma\psi&=-\gamma\varphi+\frac{n+m-2}{4(n+m-1)}R_{\phi_0}^m(\varphi-\varphi^{\frac{n+m+2}{n+m-2}})~~\mbox{ in }M,\\
\frac{\partial\psi}{\partial\nu_{g_0}}-\rho\psi&=-\rho\varphi+\frac{n+m-2}{2(n+m-1)}H_{\phi_0}^m(\varphi^{\frac{n+m}{n+m-2}}-\varphi)~~\mbox{ on }\partial M,
\end{split}
\end{equation*}
where
\begin{equation}\label{2.10}
\gamma\geq\frac{n+m}{2(n+m-1)}\|R_{\phi_0}^m\|_\infty~~\mbox{ and }~~
\rho\leq\frac{1}{n+m-1}\inf_{\partial M}H_{\phi_0}^m.
\end{equation}

Elliptic regularity theory (c.f.
\cite{ADN,S}) guarantees that the map $T$ is well defined and is
compact. Assume that $H_{\phi_0}^m\leq 0$ and consider the set
$A=\{u_0\leq u\leq w\}$, where $u_0$ and $w$ are
respectively the lower and upper solutions constructed in Lemma \ref{lem2.6}.
Recall that $u_0$ depends on $\epsilon$ and $\alpha$,
while $w$ depends on $\delta$. Thus, $A=A(\epsilon,\alpha,\delta)$ is
a set depending on $\epsilon$, $\alpha$, and $\delta$.

\begin{lemma}\label{lem2.7}
The map $T$ satisfies $T(A)\subset int(A)$.
\end{lemma}
\begin{proof}
First, we show that if $u,v\in A$ with $u\geq v$,
then $Tu\geq Tv$.
To see this, since $1\geq u\geq v\geq 0$,
for any $a\geq 0$ we have
$$u^a-v^a\leq a(u-v).$$
Using this inequality, we have in $M$
\begin{equation}\label{2.11}
\begin{split}
&\left[\gamma u-\frac{n+m-2}{4(n+m-1)}R_{\phi_0}^m(u-u^{\frac{n+m+2}{n+m-2}})\right]
-\left[\gamma v-\frac{n+m-2}{4(n+m-1)}R_{\phi_0}^m(v-v^{\frac{n+m+2}{n+m-2}})\right]\\
&\geq (u-v)\left(\gamma-\frac{n+m}{2(n+m-1)}\|R_{\phi_0}^m\|_\infty\right),
\end{split}
\end{equation}
and, since $H_{\phi_0}^m\leq 0$, we have on $\partial M$
\begin{equation}\label{2.12}
\begin{split}
&\left[\rho u-\frac{n+m-2}{2(n+m-1)}H_{\phi_0}^m(u^{\frac{n+m}{n+m-2}}-u)\right]
-\left[\rho v-\frac{n+m-2}{2(n+m-1)}H_{\phi_0}^m(v^{\frac{n+m}{n+m-2}}-v)\right]\\
&\leq (u-v)\left(\rho-\frac{1}{n+m-1}H_{\phi_0}^m\right).
\end{split}
\end{equation}
Thus, if $Tu=\psi_1$ and $Tv=\psi_2$, we get that $\psi_1$ solves the problem
\begin{equation*}
\begin{split}
\Delta_{\phi_0}\psi_1-\gamma\psi_1&=-\gamma u+\frac{n+m-2}{4(n+m-1)}R_{\phi_0}^m(u-u^{\frac{n+m+2}{n+m-2}})~~\mbox{ in }M,\\
\frac{\partial\psi_1}{\partial\nu_{g_0}}-\rho\psi_1&=-\rho u+\frac{n+m-2}{2(n+m-1)}H_{\phi_0}^m(u^{\frac{n+m}{n+m-2}}-u)~~\mbox{ on }\partial M,
\end{split}
\end{equation*}
while $\psi_2$ solves the problem
\begin{equation*}
\begin{split}
\Delta_{\phi_0}\psi_2-\gamma\psi_2&=-\gamma v+\frac{n+m-2}{4(n+m-1)}R_{\phi_0}^m(v-v^{\frac{n+m+2}{n+m-2}})~~\mbox{ in }M,\\
\frac{\partial\psi_2}{\partial\nu_{g_0}}-\rho\psi_2&=-\rho v+\frac{n+m-2}{2(n+m-1)}H_{\phi_0}^m(v^{\frac{n+m}{n+m-2}}-v)~~\mbox{ on }\partial M.
\end{split}
\end{equation*}
Therefore, it follows from (\ref{2.10})-(\ref{2.12}) that
the function $\psi=\psi_2-\psi_1$ satisfies
\begin{equation}
\begin{split}
\Delta_{\phi_0}\psi-\gamma\psi&\geq 0~~\mbox{ in }M,\\
\frac{\partial\psi}{\partial\nu_{g_0}}-\rho\psi&\leq 0~~\mbox{ on }\partial M.
\end{split}
\end{equation}
The maximum principle then implies that $\psi\leq 0$, or equivalently,
$Tv=\psi_2\leq\psi_1=T u$, as required.

Note that $u_0$ and $w$ are respectively strict lower and upper solutions.
Hence, it follows from the maximum principle that
$Tu_0> u_0$ and $Tw<w$. Therefore, if $v\in A$, then
$$u_0<T u_0\leq Tv\leq Tw<w,$$
which proves the assertion.
\end{proof}

A consequence of Lemma \ref{lem2.7} is that $T$ has no fixed point on $\partial A$, the boundary of $A$.
This allows us to define the Leray-Schauder degree of the map $I-T$ on the set $A$, which will
be denoted by $\deg(I-T,A,0)$.

\begin{lemma}\label{lem2.8}
$\deg(I-T,A,0)=1$.
\end{lemma}
\begin{proof}
Let $v_0\in int(A)$ and consider the consider map $C_{v_0}: A\to A$ defined by
$C_{v_0}(u)=v_0$. We define
$F(t,u)=t Tu+(1-t) C_{v_0}(u)$, where $F:[0,1]\times A\to A$. We claim that $F$ has no fixed point
on $\partial A$. To see this, let $v_1\in A$ and $t_0$ with $0<t_0\leq 1$ such that
$F(t_0,v_1)=v_1$, or equivalently,
\begin{equation}\label{2.13}
t_0 Tv_1+(1-t_0)v_0=v_1.
\end{equation}
Since $v_1\leq w$, we have $Tv_1\leq Tw<w$.
Then (\ref{2.13}) implies that $v_1<w$.
Similarly, since $v_1\geq u_0$, we have $Tv_1\geq T u_0>u_0$.
Then (\ref{2.13}) implies that $v_1>u_0$.
Thus, $v_1\in int(A)$.

The homotopy invariance of the Leray-Schauder degree asserts that
$\deg(I-F(t,\cdot), A,0)$ is independent of $t$. Therefore, we have
$$\deg(I-T, A,0)=\deg(I-F(1,\cdot), A,0)=\deg(I-F(0,\cdot), A,0)
=\deg(I-C_{v_0}, A,0)=1,$$
as required.
\end{proof}

\begin{lemma}\label{lem2.9} If $u$ is a solution to \eqref{2.1},
then $u$ is an isolated fixed point of $T$ and the index of $T$ at $u$,
which is defined as
$$i(T,u)=\deg(I-T,B_\delta(u),0),$$
where $B_\delta(u)=\{\varphi\in C^\infty(M):\|\varphi-u\|_\infty<\delta\}$
and $\delta$ is a small positive real number, is equal to $1$.
\end{lemma}
\begin{proof}
First we claim that the derivative of $T$ at $u$, $D_uT$, has no eigenvalue $\lambda\geq 1$.
Observe that $(D_uT)\varphi$ is defined as the unique solution to the problem:
\begin{equation*}
\begin{split}
(\Delta_{\phi_0}-\gamma)\big((D_uT)\varphi\big)&=\left(-\gamma+\frac{n+m-2}{4(n+m-1)}R_{\phi_0}^m\left(1-\frac{n+m+2}{n+m-2}u^{\frac{4}{n+m-2}}\right)\right)\varphi~~\mbox{ in }M,\\
\left(\frac{\partial}{\partial\nu_{g_0}}-\rho\right)\big((D_uT)\varphi\big)&=\left(-\rho+\frac{n+m-2}{2(n+m-1)}H_{\phi_0}^m\left(\frac{n+m}{n+m-2}u^{\frac{2}{n+m-2}}-1\right)\right)\varphi~~\mbox{ on }\partial M,
\end{split}
\end{equation*}
In particular, if $(D_uT)\varphi=\lambda\varphi$, then $\varphi$ satisfies
\begin{equation*}
\begin{split}
\lambda\Delta_{\phi_0}\varphi&=\left(\gamma(\lambda-1)+\frac{n+m-2}{4(n+m-1)}R_{\phi_0}^m\left(1-\frac{n+m+2}{n+m-2}u^{\frac{4}{n+m-2}}\right)\right)\varphi~~\mbox{ in }M,\\
\lambda\frac{\partial\varphi}{\partial\nu_{g_0}}&=\left(\rho(\lambda-1)+\frac{n+m-2}{2(n+m-1)}H_{\phi_0}^m\left(\frac{n+m}{n+m-2}u^{\frac{2}{n+m-2}}-1\right)\right)\varphi~~\mbox{ on }\partial M,
\end{split}
\end{equation*}
Multiplying the first equation by $\varphi$ and integrating by parts, we find
\begin{equation*}
\begin{split}
0&=\lambda\int_M|\nabla_{g_0}\varphi|^2e^{-\phi_0}dV_{g_0}
+(\lambda-1)\gamma\int_M\varphi^2e^{-\phi_0}dV_{g_0}-(\lambda-1)\rho\int_{\partial M}\varphi^2e^{-\phi_0}dA_{g_0}\\
&\hspace{4mm}+\frac{n+m-2}{4(n+m-1)}\int_MR_{\phi_0}^m\left(1-\frac{n+m+2}{n+m-2}u^{\frac{4}{n+m-2}}\right)\varphi^2e^{-\phi_0}dV_{g_0}\\
&\hspace{4mm}-\frac{n+m-2}{2(n+m-1)}\int_{\partial M}H_{\phi_0}^m\left(\frac{n+m}{n+m-2}u^{\frac{2}{n+m-2}}-1\right)e^{-\phi_0}dA_{g_0}.
\end{split}
\end{equation*}
Therefore, if $\lambda\geq 1$, by \eqref{2.10} we have
\begin{equation*}
\begin{split}
0&\geq \int_M|\nabla_{g_0}\varphi|^2e^{-\phi_0}dV_{g_0}+\frac{n+m-2}{4(n+m-1)}\int_MR_{\phi_0}^m\left(1-\frac{n+m+2}{n+m-2}u^{\frac{4}{n+m-2}}\right)\varphi^2e^{-\phi_0}dV_{g_0}\\
&\hspace{4mm}-\frac{n+m-2}{2(n+m-1)}\int_{\partial M}H_{\phi_0}^m\left(\frac{n+m}{n+m-2}u^{\frac{2}{n+m-2}}-1\right)e^{-\phi_0}dA_{g_0}.
\end{split}
\end{equation*}
Hence, there exists an eigenfunction $\varphi_1>0$ and an eigenvalue $\lambda_1<0$
satisfying
\begin{equation}\label{2.14}
\begin{split}
\Delta_{\phi_0}\varphi_1
-\frac{n+m-2}{4(n+m-1)}R_{\phi_0}^m\left(1-\frac{n+m+2}{n+m-2}u^{\frac{4}{n+m-2}}\right)\varphi_1+\lambda_1\varphi_1&=0~~\mbox{ in }M,\\
\frac{\partial\varphi_1}{\partial\nu_{g_0}}-\frac{n+m-2}{2(n+m-1)}H_{\phi_0}^m\left(\frac{n+m}{n+m-2}u^{\frac{2}{n+m-2}}-1\right)\varphi_1&=0~~\mbox{ on }\partial M.
\end{split}
\end{equation}
Multiplying the first equation of (\ref{2.14}) by $\big(u-u^{\frac{n+m+2}{n+m-2}}\big)>0$, we obtain
\begin{equation*}
\begin{split}
0&=\int_M\big(u-u^{\frac{n+m+2}{n+m-2}}\big)\Delta_{\phi_0}\varphi_1 e^{-\phi_0}dV_{g_0}+\lambda_1\int_M\varphi_1\big(u-u^{\frac{n+m+2}{n+m-2}}\big)e^{-\phi_0}dV_{g_0}\\
&\hspace{4mm}
-\frac{n+m-2}{4(n+m-1)}\int_MR_{\phi_0}^m\left(1-\frac{n+m+2}{n+m-2}u^{\frac{4}{n+m-2}}\right)\varphi_1\big(u-u^{\frac{n+m+2}{n+m-2}}\big)e^{-\phi_0}dV_{g_0}.
\end{split}
\end{equation*}
Using the first equation in (\ref{2.1}), it can be written as
\begin{equation}\label{2.15}
\begin{split}
0&=\int_M\big(u-u^{\frac{n+m+2}{n+m-2}}\big)\Delta_{\phi_0}\varphi_1 e^{-\phi_0}dV_{g_0}+\lambda_1\int_M\varphi_1\big(u-u^{\frac{n+m+2}{n+m-2}}\big)e^{-\phi_0}dV_{g_0}\\
&\hspace{4mm}
-\int_M\left(1-\frac{n+m+2}{n+m-2}u^{\frac{4}{n+m-2}}\right)\varphi_1\Delta_{\phi_0}u e^{-\phi_0}dV_{g_0}.
\end{split}
\end{equation}
Integrating by parts and using the second equations in (\ref{2.1}) and in (\ref{2.14}),
we obtain
\begin{equation*}
\begin{split}
&\int_M u\Delta_{\phi_0}\varphi_1 e^{-\phi_0}dV_{g_0}-\int_M \varphi_1\Delta_{\phi_0}u e^{-\phi_0}dV_{g_0}\\
&=\int_{\partial M}u\frac{\partial\varphi_1}{\partial\nu_{g_0}}e^{-\phi_0}dA_{g_0}-\int_{\partial M}\varphi_1\frac{\partial u}{\partial\nu_{g_0}}e^{-\phi_0}dA_{g_0}
=\frac{1}{n+m-1}\int_{\partial M}H_{\phi_0}^m u^{\frac{n+m}{n+m-2}}\varphi_1 e^{-\phi_0}dA_{g_0}.
\end{split}
\end{equation*}
From this and (\ref{2.15}), we get
\begin{equation*}
\begin{split}
&-\int_Mu^{\frac{n+m+2}{n+m-2}}\Delta_{\phi_0}\varphi_1 e^{-\phi_0}dV_{g_0}+
\frac{n+m+2}{n+m-2}\int_Mu^{\frac{4}{n+m-2}}\varphi_1\Delta_{\phi_0}u e^{-\phi_0}dV_{g_0}\\
&\geq -\frac{1}{n+m-1}\int_{\partial M}H_{\phi_0}^m u^{\frac{n+m}{n+m-2}}\varphi_1 e^{-\phi_0}dA_{g_0}.
\end{split}
\end{equation*}
Integrating by parts yields
\begin{equation*}
\begin{split}
&-\frac{4(n+m+2)}{(n+m-2)^2}\int_Mu^{\frac{4}{n+m-2}-1}|\nabla_{\phi_0}u|^2\varphi_1 e^{-\phi_0}dV_{g_0}
-\int_{\partial M}u^{\frac{n+m+2}{n+m-2}}\frac{\partial\varphi_1}{\partial\nu_{g_0}} e^{-\phi_0}dA_{g_0}\\
&+\frac{n+m+2}{n+m-2}\int_{\partial M}u^{\frac{4}{n+m-2}}\varphi_1\frac{\partial u}{\partial\nu_{g_0}} e^{-\phi_0}dA_{g_0}
\\
&\hspace{4mm}\geq -\frac{1}{n+m-1}\int_{\partial M}H_{\phi_0}^m u^{\frac{n+m}{n+m-2}}\varphi_1 e^{-\phi_0}dA_{g_0}.
\end{split}
\end{equation*}
This together with the second equations in (\ref{2.1}) and in (\ref{2.15}) gives
\begin{equation*}
\begin{split}
&\frac{4(n+m+2)}{(n+m-2)^2}\int_Mu^{\frac{4}{n+m-2}-1}|\nabla_{\phi_0}u|^2\varphi_1 e^{-\phi_0}dV_{g_0}\\
&\leq \frac{1}{n+m-1}\int_{\partial M}H_{\phi_0}^m u^{\frac{n+m}{n+m-2}}\varphi_1 e^{-\phi_0}dA_{g_0}
-\int_{\partial M}u^{\frac{n+m+2}{n+m-2}}\frac{\partial\varphi_1}{\partial\nu_{g_0}} e^{-\phi_0}dA_{g_0}\\
&\hspace{4mm}+\frac{n+m+2}{n+m-2}\int_{\partial M}u^{\frac{4}{n+m-2}}\varphi_1\frac{\partial u}{\partial\nu_{g_0}} e^{-\phi_0}dA_{g_0}\\
&=\frac{1}{n+m-1}\int_{\partial M}H_{\phi_0}^m u^{\frac{n+m}{n+m-2}}\left(1-u^{\frac{2}{n+m-2}}+u^{\frac{4}{n+m-2}}\right)\varphi_1 e^{-\phi_0}dA_{g_0}\\
&\leq \frac{1}{n+m-1}\int_{\partial M}H_{\phi_0}^m u^{\frac{n+m}{n+m-2}}\left(1-u^{\frac{2}{n+m-2}}\right)^2\varphi_1 e^{-\phi_0}dA_{g_0}
\leq 0
\end{split}
\end{equation*}
where the last two inequalities follow from the fact that $\varphi_1>0$ and $H_{\phi_0}^m\leq 0$.
This implies that $u$ is constant, which is a contradiction.

Therefore, all eigenvalue of $D_uT$ are small than $1$,
and thus we have that $i(T,u)=1$ (c.f. \cite{Nirenberg}).
\end{proof}

\begin{proof}[Proof of Theorem \ref{thm3}]
If there exists $[(g,\phi)]<[(g_0,\phi_0)]$ such that $R^m_{\phi}=R^m_{\phi_0}$ in $M$ and
$H^m_{\phi}=H^m_{\phi_0}\leq 0$ on $\partial M$,
then it follows from Lemma \ref{lem2.4} and Lemma \ref{lem2.5} that
$\lambda_1(L_{\phi_0}^m, B_{\phi_0}^m)<0$ and $\lambda_1(\overline{L_{\phi_0}^m}, \overline{B_{\phi_0}^m})<0$.
Conversely, if
$\lambda_1(L_{\phi_0}^m, B_{\phi_0}^m)<0$ and $\lambda_1(\overline{L_{\phi_0}^m}, \overline{B_{\phi_0}^m})<0$,
then Lemma \ref{lem2.6} implies the existence of $[(g,\phi)]<[(g_0,\phi_0)]$
with with $R^m_{\phi}=R^m_{\phi_0}$ in $M$ and
$H^m_{\phi}=H^m_{\phi_0}\leq 0$ on $\partial M$.

To prove the uniqueness, we
consider
$$
(M, g_i, e^{-\phi_i}dV_{g_i},e^{-\phi_i}dA_{g_i}, m)=(M,u_i^{\frac{4}{m+n-2}}g_0,u_i^{\frac{2(m+n)}{m+n-2}}e^{-\phi_0}dV_{g_0},u_i^{\frac{2(m+n)}{m+n-2}}e^{-\phi_0}dA_{g_0}, m)$$
satisfying $R_{\phi_i}^m=R_{\phi_0}^m$ in $M$
and $H_{\phi_i}^m=H_{\phi_0}^m$ on $\partial M$
with $0<u_i<1$ for $i=1,2$.
Choosing $\epsilon$, $\alpha$, and $\delta$
such that $u_0<u_i<w$.
Lemma \ref{lem2.9} asserts that
any solution $u$ to \eqref{2.1}
in $A=A(\epsilon,\alpha,\delta)$ is an isolated fixed point of
$T$. Thus a standard compactness argument implies that there exists
at most finite number of solutions to problem (\ref{2.1})
in $A$. Let $u_1,..., u_k$ be those solutions.
Lemma \ref{lem2.8}, Lemma \ref{lem2.9} and
the additivity of the Leray-Schauder degree imply that
$$1=\deg(I-T,A,0)=\sum_{j=1}^k i(T,u_j)=k.$$
Therefore, $u_1=u_2$, or equivalently, $g_1=g_2$.
\end{proof}

We remark that if $(M,g,e^{-\phi}dV_{g},e^{-\phi}dA_{g}, m)$
in Theorem \ref{thm3} exists,
then we have $\lambda_1(\overline{L_{\phi}^m}, \overline{B_{\phi}^m})>0$.
To see this, observe that Theorem \ref{thm3} asserts that
$\lambda_1(\overline{L_{\phi_0}^m}, \overline{B_{\phi_0}^m})<0$
and Theorem \ref{thm2.2} implies that $\lambda_1(\overline{L_{\phi}^m}, \overline{B_{\phi}^m})\neq 0$.  Since the sign of
$\lambda_1(L_{\phi_0}^m, B_{\phi_0}^m)$ is conformal invariant,
Theorem \ref{thm3} implies that
$\lambda_1(L_{\phi}^m, B_{\phi}^m)<0$.
Therefore, if $\lambda_1(\overline{L_{\phi}^m}, \overline{B_{\phi}^m})<0$,
then Theorem \ref{thm3} implies that
there exists
$(M,\tilde{g},e^{-\tilde{\phi}}dV_{\tilde{g}},e^{-\tilde{\phi}}dA_{\tilde{g}}, m)$
conformal to $(M,g,e^{-\phi}dV_{g},e^{-\phi}dA_{g}, m)$
such that $\tilde{g}<g$ with $R^m_{\tilde{\phi}}=R^m_{\phi}$ in $M$ and
$H^m_{\tilde{\phi}}=H^m_{\phi}\leq 0$ on $\partial M$.
But this contradicts to the uniqueness statement of Theorem \ref{thm3}.
Hence, we must have $\lambda_1(\overline{L_{\phi}^m}, \overline{B_{\phi}^m})>0$.

In the rest of this section, we give two geometric criteria
which guarantee when $\lambda_1(\overline{L_{\phi_0}^m}, \overline{B_{\phi_0}^m})$
or $\lambda_1(L_{\phi_0}^m, B_{\phi_0}^m)$
is negative.

\begin{prop}\label{prop2.11}
Let
$(M,g_0,e^{-\phi_0}dV_{g_0},e^{-\phi_0}dA_{g_0}, m)$
be
a compact smooth metric measure space with boundary
such that $R_{\phi_0}^m$
and $H_{\phi_0}^m$
do not vanish simultaneously.
If $\displaystyle\int_MR_{\phi_0}^m e^{-\phi_0}dV_{g_0}
+\int_{\partial M}H_{\phi_0}^me^{-\phi_0}dA_{g_0}\geq 0$,
then $\lambda_1(\overline{L_{\phi_0}^m}, \overline{B_{\phi_0}^m})<0$.
\end{prop}
\begin{proof}
Let $f_1$ be the first eigenfunction corresponding to the eigenvalue
$\lambda_1(\overline{L_{\phi_0}^m}, \overline{B_{\phi_0}^m})$.
Thus $f_1$ satisfies the eigenvalue problem (\ref{2.4}).
Multiplying the first equation by $f_1^{-1}$ and integrating by parts
give
\begin{equation*}
\begin{split}
&\frac{1}{n+m-1}\int_MR_{\phi_0}^m e^{-\phi_0}dV_{g_0}\\
&=-\int_M\frac{\Delta_{\phi_0}f_1}{f_1}e^{-\phi_0}dV_{g_0}
-\lambda_1(\overline{L_{\phi_0}^m}, \overline{B_{\phi_0}^m})\int_Me^{-\phi_0}dV_{g_0}\\
&=-\lambda_1(\overline{L_{\phi_0}^m}, \overline{B_{\phi_0}^m})\int_Me^{-\phi_0}dV_{g_0}
-\int_M\frac{|\nabla_{g_0}f_1|^2}{f_1^2}e^{-\phi_0}dV_{g_0}
-\int_{\partial M}\frac{1}{f_1}\frac{\partial f_1}{\partial\nu_{g_0}}e^{-\phi_0}dA_{g_0}\\
&=-\lambda_1(\overline{L_{\phi_0}^m}, \overline{B_{\phi_0}^m})\int_Me^{-\phi_0}dV_{g_0}
-\int_M\frac{|\nabla_{g_0}f_1|^2}{f_1^2}e^{-\phi_0}dV_{g_0}
-\frac{1}{n+m-1}\int_{\partial M}H_{\phi_0}^me^{-\phi_0}dA_{g_0}.
\end{split}
\end{equation*}
Hence, if $\lambda_1(\overline{L_{\phi_0}^m}, \overline{B_{\phi_0}^m})\geq 0$,
we have
$$0\leq\frac{1}{n+m-1}
\left(\int_MR_{\phi_0}^m e^{-\phi_0}dV_{g_0}
+\int_{\partial M}H_{\phi_0}^me^{-\phi_0}dA_{g_0}\right)
\leq -\int_M\frac{|\nabla_{g_0}f_1|^2}{f_1^2}e^{-\phi_0}dV_{g_0}.$$
This implies that $f_1$ is a constant function.
Since $f_1$ is a constant, it follows from (\ref{2.4}) that
$H_{\phi_0}^m\equiv 0$ and
$R_{\phi_0}^m\equiv (n+m-1)\lambda_1(\overline{L_{\phi_0}^m}, \overline{B_{\phi_0}^m})\leq 0$.
Combining this with the assumption that
$\displaystyle\int_MR_{\phi_0}^m e^{-\phi_0}dV_{g_0}
+\int_{\partial M}H_{\phi_0}^me^{-\phi_0}dA_{g_0}\geq 0$,
we have $R_{\phi_0}^m\equiv 0$,
which contradicts to the assumption that  $R_{\phi_0}^m$
and $H_{\phi_0}^m$
do not vanish simultaneously.
Therefore, we must have $\lambda_1(\overline{L_{\phi_0}^m}, \overline{B_{\phi_0}^m})<0$.
\end{proof}

\begin{prop}\label{prop2.12}
Let
$(M,g_0,e^{-\phi_0}dV_{g_0},e^{-\phi_0}dA_{g_0}, m)$
be
a compact smooth metric measure space with boundary
such that $R_{\phi_0}^m$
and $H_{\phi_0}^m$
do not vanish simultaneously.
If $\displaystyle\int_MR_{\phi_0}^m e^{-\phi_0}dV_{g_0}
+2\int_{\partial M}H_{\phi_0}^me^{-\phi_0}dA_{g_0}\leq 0$,
then $\lambda_1(L_{\phi_0}^m, B_{\phi_0}^m)<0$.
\end{prop}
\begin{proof}
Let $\varphi_1$ be the first eigenfunction
of $(L_{\phi_0}^m, B_{\phi_0}^m)$.
Then $\varphi_1$ satisfies the boundary value problem
(\ref{2.9}).
Multiplying the first equation by $\varphi_1^{-1}$
and integrating by parts yield
\begin{equation*}
\begin{split}
&\int_MR_{\phi_0}^m e^{-\phi_0}dV_{g_0}\\
&=\frac{4(n+m-1)}{n+m-2}\int_M\frac{\Delta_{\phi_0}\varphi_1}{\varphi_1}e^{-\phi_0}dV_{g_0}
+\lambda_1(L_{\phi_0}^m, B_{\phi_0}^m)\int_Me^{-\phi_0}dV_{g_0}\\
&=\frac{4(n+m-1)}{n+m-2}\int_M\frac{|\nabla_{\phi_0}\varphi_1|^2}{\varphi_1^2}e^{-\phi_0}dV_{g_0}
+\frac{4(n+m-1)}{n+m-2}\int_{\partial M}\frac{1}{\varphi_1}\frac{\partial \varphi_1}{\partial\nu_{g_0}}e^{-\phi_0}dA_{g_0}\\
&\hspace{4mm}
+\lambda_1(L_{\phi_0}^m, B_{\phi_0}^m)\int_Me^{-\phi_0}dV_{g_0}\\
&=\frac{4(n+m-1)}{n+m-2}\int_M\frac{|\nabla_{\phi_0}\varphi_1|^2}{\varphi_1^2}e^{-\phi_0}dV_{g_0}
-2\int_{\partial M}H_{\phi_0}^m e^{-\phi_0}dA_{g_0}\\
&\hspace{4mm}
+\lambda_1(L_{\phi_0}^m, B_{\phi_0}^m)\int_Me^{-\phi_0}dV_{g_0}.
\end{split}
\end{equation*}
Hence, if $\lambda_1(L_{\phi_0}^m, B_{\phi_0}^m)\geq 0$,
then we have
$$\frac{4(n+m-1)}{n+m-2}\int_M\frac{|\nabla_{\phi_0}\varphi_1|^2}{\varphi_1^2}e^{-\phi_0}dV_{g_0}
\leq
\int_MR_{\phi_0}^m e^{-\phi_0}dV_{g_0}
+2\int_{\partial M}H_{\phi_0}^m e^{-\phi_0}dA_{g_0}\leq 0,$$
where we have used the assumption in the last inequality.
This implies that $\varphi_1$ is constant,
and we get a contradiction as in  the proof of Proposition \ref{prop2.11}.
Therefore, we have $\lambda_1(L_{\phi_0}^m, B_{\phi_0}^m)<0.$
\end{proof}

Combining Propositions \ref{prop2.11}, \ref{prop2.12}
and Theorem \ref{thm3}, we obtain the following:

\begin{cor}
Suppose that $(M,g_0,e^{-\phi_0}dV_{g_0},e^{-\phi_0}dA_{g_0}, m)$  is a SMMS with boundary
such that $H_{\phi_0}^m\leq 0$. There exists $(M,g,e^{-\phi}dV_{g},e^{-\phi}dA_{g}, m)$ conformal to
$(M,g_0,e^{-\phi_0}dV_{g_0},e^{-\phi_0}dA_{g_0}, m)$ in the sense of \eqref{1.5}
such that $[(g,\phi)]<[(g_0,\phi_0)]$ with $R^m_{\phi}=R^m_{\phi_0}$ in $M$ and
$H^m_{\phi}=H^m_{\phi_0}\leq 0$ on $\partial M$
if and only if
one of the following alternatives occurs:\\
(i) $\lambda_1(L_{\phi_0}^m, B_{\phi_0}^m)<0$
and $\displaystyle\int_MR_{\phi_0}^m e^{-\phi_0}dV_{g_0}
+\int_{\partial M}H_{\phi_0}^me^{-\phi_0}dA_{g_0}>0$,\\
(ii)  $\displaystyle\int_MR_{\phi_0}^m e^{-\phi_0}dV_{g_0}
+\int_{\partial M}H_{\phi_0}^me^{-\phi_0}dA_{g_0}=0$,\\
(iii) $\lambda_1(\overline{L_{\phi_0}^m}, \overline{B_{\phi_0}^m})<0$
and $\displaystyle\int_MR_{\phi_0}^m e^{-\phi_0}dV_{g_0}
+\int_{\partial M}H_{\phi_0}^me^{-\phi_0}dA_{g_0}<0$,\\
(iv) $\displaystyle\int_MR_{\phi_0}^m e^{-\phi_0}dV_{g_0}
+\int_{\partial M}H_{\phi_0}^me^{-\phi_0}dA_{g_0}=0$.
\end{cor}

\section{weighted Yamabe solitons with boundary}\label{section7}

In this section, we give the characterization of the weighted Yamabe solitons with boundary from two points of view. 

Firstly, note that the unnormalized weighted Yamabe flow with boudnary is equivalent to the (normalized) weighted Yamabe flow with boundary defined as
\begin{equation}\label{flow}
\left\{
\begin{array}{ll}
\frac{\partial g}{\partial t} &=(r^m_{\phi}-R^m_{\phi})g \\
\frac{\partial \phi}{\partial t} &=\frac{m}{2}(R^m_{\phi}-r^m_{\phi})
\end{array}
\right. ~~\mbox{ in }M ~~\mbox{ and }~~
H^m_{\phi}=0 ~~\mbox{ on }\partial M,
\end{equation}
where $r_{\phi}^m$ is the average of the weighted scalar curvature $R_{\phi}^m$; i.e.
\begin{equation}\label{3.3}
r_{\phi}^m=\frac{\int_M R_{\phi}^m e^{-\phi}dV_{g}}
{\int_M  e^{-\phi}dV_{g}}.
\end{equation}
To see this, we let $(M,g(t),e^{-\phi(t)}dV_{g(t)},e^{-\phi(t)}dA_{g(t)}, m)$ be a solution to the unnormalized weighted Yamabe flow with boundary (\ref{unnormalized}).
Suppose $\psi(t)$ is a positive function defined by 
\begin{equation}\label{4.3}
\psi'(t)=r_{\phi(t)}^m \psi(t)~~\mbox{ with }\psi(0)=1,
\end{equation}
where $r_{\phi(t)}^m$ is the average of the weighted scalar curvature $R_{\phi(t)}^m$:
\begin{equation*}
r_{\phi(t)}^m=\frac{\int_MR_{\phi(t)}^m e^{-\phi(t)}dV_{g(t)}}{\int_Me^{-\phi(t)}dV_{g(t)}}.
\end{equation*}
Let
\begin{equation}\label{4.5}
\tilde{t}=\int_0^t\psi(\tau)d\tau.
\end{equation}
Then $\displaystyle\frac{d\tilde{t}}{dt}=\psi(t)>0$, which implies that
$t$ is a differentiable function in $\tilde{t}$ by the inverse function theorem.
We claim that
 $(M,\tilde{g}(\tilde{t}),e^{-\tilde{\phi}(\tilde{t})}dV_{\tilde{g}(\tilde{t})}, e^{-\tilde{\phi}(\tilde{t})}dA_{g(\tilde{t})}, m)$ is solution of the (normalized) weighted Yamabe flow with boundary,
where
\begin{equation}\label{4.6}
\tilde{g}(\tilde{t})=\psi(t)g(t)~~\mbox{ and }~~
\tilde{\phi}(\tilde{t})=\phi(t)-\frac{m}{2}\log\psi(t).
\end{equation}
To prove this claim, we first observe that
\begin{equation}\label{4.7}
R_{\tilde{\phi}(\tilde{t})}^m=\psi(t)^{-1}R_{\phi(t)}^m~~\mbox{ and }
~~r_{\tilde{\phi}(\tilde{t})}^m=\psi(t)^{-1}r_{\phi(t)}^m
\end{equation}
by (\ref{4.6}).
We then compute
\begin{equation}\label{4.8}
\begin{split}
\frac{\partial}{\partial\tilde{t}}\tilde{g}(\tilde{t})
&=\frac{\partial}{\partial t}\left(\psi(t)g(t)\right)\frac{dt}{d\tilde{t}}\\
&=\left(\psi'(t)g(t)+\psi(t)\frac{\partial}{\partial t}g(t)\right)\frac{dt}{d\tilde{t}}\\
&=\frac{1}{\psi(t)}\left(r_{\phi(t)}^m\psi(t)g(t)-\psi(t)R_{\phi(t)}^mg(t)\right)\\
&=\left(r_{\phi(t)}^m-R_{\phi(t)}^m\right)g(t)\\
&=\left(r_{\tilde{\phi}(\tilde{t})}^m-R_{\tilde{\phi}(\tilde{t})}^m\right)\tilde{g}(\tilde{t}),
\end{split}
\end{equation}
where we have used (\ref{unnormalized})-(\ref{4.5}) and (\ref{4.7}). Similarly, it follows from
(\ref{unnormalized})-(\ref{4.5}) and (\ref{4.7}) that
\begin{equation}\label{4.9}
\begin{split}
\frac{\partial}{\partial\tilde{t}}\tilde{\phi}(\tilde{t})
&=\frac{\partial}{\partial t}\left(\phi(t)-\frac{m}{2}\log\psi(t)\right)\frac{dt}{d\tilde{t}}\\
&=\left(\frac{\partial}{\partial t}\phi(t)-\frac{m}{2}\frac{\psi'(t)}{\psi(t)}\right)\frac{dt}{d\tilde{t}}\\
&=\frac{1}{\psi(t)}\left(\frac{m}{2}R_{\phi(t)}^m-\frac{m}{2}r_{\phi(t)}^m \right)\\
&=\frac{m}{2}\left(R_{\tilde{\phi}(\tilde{t})}^m-r_{\tilde{\phi}(\tilde{t})}^m\right).
\end{split}
\end{equation}
Moreover, we have
\begin{equation}\label{eq4.10}
  H_{\tilde{\phi}(\tilde{t})}^m=H_{\tilde{g}}+\frac{\partial \tilde{\phi}}{\partial \nu_{\tilde{g}}}=\psi^{-\frac{1}{2}}\left(H_g+\frac{\partial \phi}{\partial \nu_g}\right)=\psi^{-\frac{1}{2}}H_\phi^m=0.
\end{equation}
The equivalence follows from (\ref{4.8}), (\ref{4.9}) and \eqref{eq4.10} immediately.

\medskip

Since the weighted Yamabe flow preserves the conformal structure, we may write \begin{equation}
\begin{cases}
&g(t)=w(t)^{\frac{4}{n+m-2}}g_0, \\
&e^{-\phi(t)}=w(t)^{\frac{2m}{n+m-2}}e^{-\phi_0},
\end{cases}
\end{equation} 
as the solution of (\ref{unnormalized}) with $(g(0), \phi(0))=(g_0, \phi_0)$. Hence the conformal factor $w(t)$ has the evolution equation \begin{equation}\label{sol5.4}
  \left\{
    \begin{array}{ll}
      \frac{\partial w(t)}{\partial t}=-\frac{m+n-2}{4}R_{\phi(t)}^mw(t), & \hbox{in $M$;} \\
      \frac{\partial w(t)}{\partial \nu_{g_0}}=0, & \hbox{on $\partial M$;} \\
      w(0)=1, & \hbox{}
    \end{array}
  \right.
\end{equation}
by $H_{\phi_0}^m=0$ and (\ref{1.6}).

The characterization is given in the following.

\begin{theorem}\label{thm7.1}
Any compact weighted Yamabe soliton with boundary \eqref{eq2.8} must have constant weighted scalar curvature in $M$ and vanishing weighted mean curvature on $\partial M$.
\end{theorem}
\begin{proof}
On $(M, g, e^{-\phi}dV_g,e^{-\phi}dA_g, m)$, we define the normalized energy functional $\tilde{E}(w)$ as
\begin{equation}\label{energy}
  \tilde{E}(w)=\frac{\int_MR_\phi^me^{-\phi}dV_g+2\int_{\partial M}H_\phi^me^{-\phi}dA_g}{\left(\int_M e^{-\phi}dV_{g(t)}\right)^{\frac{n+m-2}{n+m}}}.
\end{equation}
and the unnormalized version as
\begin{equation}
  E(w)=\int_MR_\phi^me^{-\phi}dV_g+2\int_{\partial M}H_\phi^me^{-\phi}dA_g.
\end{equation}
Along the flow (\ref{unnormalized}), it follows from \eqref{1.6}, \eqref{sol5.4}, and integration by parts that

\begin{equation}\label{sol5.7}
  \begin{split}
  E(w)=&\int_M R_\phi^me^{-\phi}dV_g\\
  =&\int_M R_\phi^mw^\frac{2(m+n)}{m+n-2}e^{-\phi_0}dV_{g_0}\\
  =&\int_M\left(-\frac{4(m+n-1)}{m+n-2}\Delta_{\phi_0}w+R_{\phi_0}^mw\right)we^{-\phi_0}dV_{g_0}\\
  =&\int_M \left(\frac{4(m+n-1)}{m+n-2}|\nabla w|_{g_0}^2+R_{\phi_0}^mw^2\right)e^{-\phi_0}dV_{g_0}.
\end{split}
\end{equation}

Then, by \eqref{1.6}, \eqref{sol5.4}, and integration by parts, along the weighted Yamabe flow with boundary (\ref{unnormalized}), we have

\begin{equation}\label{sol5.8}
  \begin{split}
  \frac{d}{dt}E(w)=&2\int_M\left(\frac{4(m+n-1)}{m+n-2}\langle\nabla w,\nabla w_t\rangle_{g_0}^2+R_{\phi_0}^mww_t\right)e^{-\phi_0}dV_{g_0}\\
=&\int_M\left(-\frac{4(m+n-1)}{m+n-2}\Delta_{\phi_0}w+R_{\phi_0}^mw\right)w_te^{-\phi_0}dV_{g_0}\\
&\quad +\frac{4(m+n-1)}{m+n-2}\int_{\partial M}w_t\frac{\partial w}{\partial \nu_{g_0}}e^{-\phi_0}dA_{g_0}\\
=&-\frac{m+n-2}{4}\int_M (R_\phi^m)^2w^\frac{2(m+n)}{m+n-2}e^{-\phi_0}dV_{g_0}\\
=&-\frac{m+n-2}{2}\int_M (R_\phi^m)^2 e^{-\phi}dV_g.
\end{split}
\end{equation}

Note also that

\begin{equation}\label{sol5.9}
  \frac{d}{dt}\left(\int_M e^{-\phi}dV_g\right)=\frac{d}{dt}\left(\int_M w^\frac{2(m+n)}{m+n-2}e^{-\phi_0}dV_{g_0}\right)=-\frac{m+n}{2}\int_M R_\phi^me^{-\phi}dV_g
\end{equation}
where the last equality follows from \eqref{sol5.4}. Combining \eqref{sol5.7}, \eqref{sol5.8}, and \eqref{sol5.9}, we obtain
\begin{equation}\label{sol5.10}
\begin{split}
&\frac{d}{dt}(\tilde{E}(w(t)))=\frac{d}{dt}\left(\frac{\int_M R_{\phi(t)}^m e^{-\phi(t)}dV_{g(t)}}{\left(\int_M e^{-\phi(t)}dV_{g(t)}\right)^{\frac{n+m-2}{n+m}}}\right)\\
&=-\frac{n+m-2}{2}\frac{\left(\int_M (R_{\phi(t)}^m)^2e^{-\phi(t)}dV_{g(t)}\right)\left(\int_M e^{-\phi(t)}dV_{g(t)}\right)-\left(\int_M R_{\phi(t)}^me^{-\phi(t)}dV_{g(t)}\right)^2}{\left(\int_M e^{-\phi(t)}dV_{g(t)}\right)^{\frac{2(n+m-1)}{n+m}}}\leq 0
\end{split}
\end{equation}
where the last inequality follows from Cauchy-Schwartz inequality.

\medskip

On the other hand, the normalized energy \eqref{energy} is invariant along the weighted Yamabe soliton \eqref{eq2.8}. To see this, note that
\begin{equation*}
 \begin{split}
       & R_{\phi(t)}^m=\sigma(t)^{-1}\psi_t^*(R_{\phi_0}^m) ~~\mbox{and}~~ e^{-\phi(t)}dV_{g(t)}=\sigma(t)^{\frac{n+m}{2}}\psi_t^*(e^{-\phi_0}dV_{g_0}) ~~\mbox{in $M$}, \\
       &H_{\phi(t)}^m=\sigma(t)^{-\frac{1}{2}}\psi_t^*(H_{\phi_0}^m) ~~\mbox{and}~~ e^{-\phi(t)}dA_{g(t)}=\sigma(t)^{\frac{n+m-1}{2}}\psi_t^*(e^{-\phi_0}dA_{g_0})~~\mbox{on $\partial M$},
   \end{split}
\end{equation*}
which implies that
\begin{equation*}
  \begin{split}
 \tilde{E}(w(t))=&\frac{\int_M R_{\phi(t)}^m e^{-\phi(t)}dV_{g(t)}+2\int_{\partial M}H_{\phi(t)}^me^{-\phi(t)}dA_{g(t)}}{\left(\int_M e^{-\phi(t)}dV_{g(t)}\right)^{\frac{n+m-2}{n+m}}}\\
 =&\frac{\sigma(t)^{\frac{n+m}{2}-1}\left(\int_M \psi_t^*(R_{\phi_0}^m)\psi_t^*(e^{-\phi_0}dV_{g_0})+2\int_{\partial M} \psi_t^*(H_{\phi_0}^m)\psi_t^*(e^{-\phi_0}dA_{g_0})\right)}{\left(\sigma(t)^{\frac{n+m}{2}}\int_M \psi_t^*(e^{-\phi_0}dV_{g_0})\right)^{\frac{n+m-2}{n+m}}}\\
=&\frac{\int_M \psi_t^*(R_{\phi_0}^me^{-\phi_0}dV_{g_0})+2\int_{\partial M} \psi_t^*(H_{\phi_0}^me^{-\phi_0}dA_{g_0})}{\left(\int_M \psi_t^*(e^{-\phi_0}dV_{g_0})\right)^{\frac{n+m-2}{n+m}}}\\
=&\frac{\int_M R_{\phi_0}^m e^{-\phi_0}dV_{g_0}+2\int_{\partial M} H_{\phi_0}^m e^{-\phi_0}dA_{g_0}}{\left(\int_M e^{-\phi_0}dV_{g_0}\right)^{\frac{n+m-2}{n+m}}}
\end{split}
\end{equation*}
where we have used the fact that $\psi_t$ is diffeomorphism of $M$.
In particular, we have
\begin{equation*}
  \frac{d}{dt} \tilde{E}(w(t))=0.
\end{equation*}
Combining this with \eqref{sol5.10}, we can conclude that
$R_{\phi(t)}^m$ is  constant by the equality case of Cauchy-Schwartz inequality. This proves the assertion.
\end{proof}

We now consider the weighted Yamabe soliton from the equation point of view.
Differentiating \eqref{eq2.8} and evaluating it at $t=0$, we obtain
\begin{equation}\label{eq2.14}
  (\lambda-R_{\phi_0}^m)g_0=\mathcal{L}_Xg_0,\quad \textrm{ and } \quad \mathcal{L}_X\phi_0=\frac{m}{2}(R_{\phi_0}^m-\lambda).
\end{equation}
where $\lambda=-\dot{\sigma}(0)$, $\mathcal{L}_X$ is the Lie derivative with respect to $X$, and $X$ is the complete vector field generated by $\psi_t$. Also, the condition that $H_{\phi}^m=0$ on $\partial M$ in \eqref{unnormalized} gives
\begin{equation}\label{sol2.10}
  H_{\phi_0}^m=0\textrm{ on } \partial M.
\end{equation}
On the other hand, the vector field $X$ satisfies the following:
\begin{equation}\label{bdcond1}
  \langle X,\nu_{g_0}\rangle =0\textrm{ on }\partial M.
\end{equation}

To see (\ref{bdcond1}), for any fixed $p\in\partial M$, consider a curve $\gamma(t):[0,T)\rightarrow \partial M$ defined by $\gamma(t):=\psi_t(p)$. Since the restriction of the $\psi_t$ to the boundary $\partial M$ is a diffeomorphism on $\partial M$, the curve $\gamma(t)$ lies on $\partial M$ for all $t\in [0,T)$. Thus we have
\begin{equation}\label{gammatnu}
  \langle\gamma'(t),\nu_{g_0}\rangle=0\quad \textrm{ on }\partial M
\end{equation}
By the construction of the vector field $X$, the left hand side of (\ref{gammatnu}) is
\begin{equation}\label{derivation}
  \langle \gamma'(t),\nu_{g_0}\rangle=\langle\frac{\partial}{\partial t}\left(\psi_t(p)\right),\nu_{g_0}\rangle=\langle X(p),\nu_{g_0}\rangle.
\end{equation}
Since $p\in\partial M$ is arbitrary, we get (\ref{bdcond1}). In view of the discussion above, we say that a  $(M,g_0,X,\phi_0,\lambda)$ is a \textit{weighted Yamabe soliton with boundary} if
\begin{equation}
  \left\{
    \begin{array}{ll}
      (\lambda-R_{\phi_0}^m)g_0=\mathcal{L}_Xg_0 &  \\
      \mathcal{L}_X\phi_0=\frac{m}{2}(R_{\phi_0}^m-\lambda) &     \end{array}
  \right.
\mbox{ in }M~~\mbox{ and }~~ H_{\phi_0}^m=0,~~X\perp \nu_{g_0}\mbox{ on }\partial M.
\end{equation}

We show in the following lemma that (\ref{2.12def}) is  actually equivalent to (\ref{eq2.8}).

\begin{lemma}\label{sollem2.1}
  If $(M,g_0,X,\phi_0,\lambda)$ satisfies \eqref{2.12def}, there exist  a solution $(M,g(t),e^{-\phi(t)}dV_{g(t)},m)$ of the unnormalized
weighted Yamabe flow
\eqref{unnormalized} satisfying
$(M,g(0),e^{-\phi(0)}dV_{g(0)},m)=(M,g_0,e^{-\phi_0}dV_{g_0},m)$,
 a family of diffeomorphisms $\psi_t$ in $M$ with $\psi_0=id_M$, a vector field $X(t)$ with $X(0)=X$ defined for all $t$ with
\begin{equation}
  \sigma(t):=1-\lambda t>0,
\end{equation}
such that the followings hold:

{\rm (i)} $\psi_t:M\to M$ is the 1-parameter family of diffeomorphisms generated by $X(t):=\frac{1}{\sigma(t)}X$,

{\rm (ii)} $g(t)$ is the pullback by $\psi_t$ of $g_0$ up to the scale factor $\sigma(t)$, and $\phi(t)=\psi_t^*(\phi_0)-\frac{m}{2}\log\sigma(t)$, and

{\rm (iii)} $X(t)$ is the pullback of $X$ by $\psi_t$.
\end{lemma}
\begin{proof}
  We define $\sigma(t)=1-\lambda t$. Since the vector field $X$ is complete and $X\perp \nu_{g_0}$, there exists a 1-parameter family of diffeomorphisms $\psi_t:M\to M$ generated by the vector fields $X(t)=\frac{1}{\sigma(t)}X$ defined for all $t$ with $\sigma(t)>0$ (see Theorem 9.34 in \cite{Lee}). Define
\begin{equation}\label{eq2.16}
  g(t)=\sigma(t) \psi_t^*(g_0),\textrm{ and }\phi(t)=\psi_t^*(\phi_0)-\frac{m}{2}\log \sigma(t).
\end{equation}
It follows from \eqref{eq2.16} that
\begin{equation}
  R_{\phi(t)}^m=\sigma(t)^{-1}\psi_t^*(R_{\phi_0}^m),
\end{equation}
which together with the first equation in  \eqref{eq2.14} implies that
\begin{equation}\label{eq2.18}
  \begin{split}
  R_{\phi(t)}^m g(t)=&\sigma(t)^{-1}\psi_t^*(R_{\phi_0}^m)\sigma(t)\psi_t^*(g_0)=\psi_t^*(R_{\phi_0}^mg_0)\\
=&\psi_t^*(\lambda g_0-\mathcal{L}_{X}g_0)=\lambda \psi_t^*(g_0)-\psi_t^*(\mathcal{L}_{X}g_0)\\
=&\frac{\lambda}{\sigma(t)}g(t)-\mathcal{L}_{X(t)}g(t).
\end{split}
\end{equation}
Differentiating the first equation in \eqref{eq2.16} with respect to $t$ and using \eqref{eq2.18}, we obtain
\begin{equation*}
  \begin{split}
  \frac{\partial}{\partial t}g(t)=&\sigma'(t)\psi_t^*(g_0)+\sigma(t)\frac{\partial}{\partial t}\psi_t^*(g_0)\\
=&-\frac{\lambda}{\sigma(t)}g(t)+\mathcal{L}_{X(t)}g(t)=-R_{\phi(t)}^mg(t).
\end{split}
\end{equation*}
By the second equation in \eqref{eq2.14}, we have
\begin{equation}\label{eq2.13}
  \begin{split}
  R_{\phi(t)}^m=&\sigma(t)^{-1}\psi_t^*(R_{\phi_0}^m)\\
=&\sigma(t)^{-1}\psi_t^*\left(\frac{2}{m}\mathcal{L}_{X}\phi_0+\lambda\right)\\
=&\frac{2}{m}\sigma(t)^{-1}\mathcal{L}_{X(t)}(\psi_t^*\phi_0)+\lambda\sigma(t)^{-1}.
\end{split}
\end{equation}
Differentiating the second  equation in \eqref{eq2.16} with respect to $t$ and
using (\ref{eq2.13}), we get
\begin{equation*}
\begin{split}
  \frac{\partial}{\partial t}\phi(t)=&\frac{\partial}{\partial t}\psi_t^*(\phi_0)-\frac{m}{2}\frac{\sigma'(t)}{\sigma(t)}\\
=&\frac{1}{\sigma(t)}\mathcal{L}_{X(t)}(\psi_t^*\phi_0)+\frac{m}{2}\frac{\lambda}{\sigma(t)}=\frac{m}{2}R_{\phi(t)}^m.
\end{split}
\end{equation*}

Also, by $H_{\phi_0}^m=0$, we have
\begin{equation*}
  H_{\phi(t)}^m=\sigma(t)^{-\frac{1}{2}}\psi_t^*(H_{\phi_0}^m)=0\textrm{ on }\partial M.
\end{equation*}

This proves the lemma.
\end{proof}

As an analogue of Theorem \ref{thm7.1}, we characterize the compact gradient weighted Yamabe soliton with boundary in the following.

\begin{theorem}\label{thm7.3}
  Any compact gradient weighted Yamabe soliton with boundary \eqref{eq2.21} with $\frac{\partial \phi_0}{\partial \nu_{g_0}}=0$ must have constant weighted scalar curvature in $M$.
\end{theorem}
\begin{proof}
First we will show that
\begin{equation}\label{sol2.22}
  \frac{\partial R_{\phi_0}^m}{\partial \nu_{g_0}}=0
\end{equation}
It follows from Lemma \ref{sollem2.1} that there exists a solution $g(t)=\sigma(t)\psi_t^*(g_0)$ and $\phi(t)=\psi_t^*(\phi_0)-\frac{m}{2}\log \sigma(t)$ to \eqref{unnormalized}. Since the flow \eqref{unnormalized} preserves conformal structure, we can write $g(t)=w(t)^\frac{4}{m+n-2}g_0$ and $e^{-\phi(t)}=w(t)^\frac{2m}{m+n-2}e^{-\phi_0}$ for some positive function $w(t)$.  Differentiate this with respect to $t$ and evaluate at $t=0$:
\begin{equation}\label{sol2.23}
  \left.\frac{\partial}{\partial t}g(t)\right|_{t=0}=\frac{4}{m+n-2}\frac{\partial w}{\partial t}(0)g_0=-R_{\phi_0}^mg_0.
\end{equation}

Since $H_{\phi(t)}^m\equiv 0$ on $\partial M$,  we have
\begin{equation*}
  H_{\phi(t)}^m=\left(\frac{2(m+n-1)}{m+n-2}\frac{\partial w(t)}{\partial \nu_{g_0}}+H_{\phi_0}^mw\right)w^\frac{m+n}{m+n-2}=0,
\end{equation*}
which implies that
\begin{equation}\label{sol2.24}
  \frac{\partial w(t)}{\partial \nu_{g_0}}=0\textrm{ for all }t.
\end{equation}
Thus we have
\begin{equation}
  \begin{split}
  \frac{\partial R_{\phi_0}^m}{\partial \nu_{g_0}}=&\frac{\partial }{\partial \nu_{g_0}}\left(-\frac{4}{m+n-2}\left.\frac{\partial w(t)}{\partial t}\right|_{t=0}\right)\\
=&-\frac{4}{m+n-2}\frac{\partial}{\partial t}\left.\left(\frac{\partial w(t)}{\partial \nu_{g_0}}\right)\right|_{t=0}=0
\end{split}
\end{equation}
where we used equalities in \eqref{sol2.23} and \eqref{sol2.24}. This proves our claim in \eqref{sol2.22}. We remark that \eqref{sol2.22} still holds for the weighted Yamabe soliton with boundary which is not gradient.

In the rest of this proof, to simplify the notation, we denote $g_0$ and $\phi_0$  as $g$ and $\phi$. And all integrals over $M$ (resp. over $\partial M$) are taken with respect to $e^{-\phi}dV_{g}$ (resp. $e^{-\phi}dA_{g}$).
Taking the trace of the first equation in \eqref{eq2.21}, we find
\begin{equation}\label{delphif}
\begin{split}
  \Delta_{\phi} f&=\Delta f-\langle \nabla f,\nabla \phi\rangle\\
&=n(\lambda-R_\phi^m)-\langle \nabla f,\nabla\phi\rangle\\
&=\left(n+m\right)(\lambda-R_{\phi}^m)
\end{split}
\end{equation}
where the last equality follows from the second equation in \eqref{eq2.21}. Integrating it over $M$ yields
\begin{equation*}
  \left(n+m\right)\int_M \left(\lambda-R_\phi^m\right)=\int_M \Delta_{\phi}f=\int_{\partial M}\frac{\partial f}{\partial \nu_{g}}=0
\end{equation*}
where the last equality follows from \eqref{eq2.21}. Since $n+m>0$, we have
\begin{equation}\label{int_trace}
  \int_M (\lambda-R_\phi^m)=0.
\end{equation}
On the other hand, we  take the covariant derivatives of the first equation in (\ref{eq2.21}) and obtain
  \begin{equation*}
    -\nabla_iR_\phi^mg_{jk}=\nabla_i\nabla_j\nabla_k f,
  \end{equation*}
which implies that
\begin{equation*}
  R_{ijkl}\nabla^l f=\nabla_i\nabla_j\nabla_k f-\nabla_j\nabla_i\nabla_k f=-\nabla_iR_\phi^mg_{jk}+\nabla_jR_\phi^mg_{ik}.
\end{equation*}
Taking the trace of $i$ and $k$ gives
\begin{equation}\label{eq2.26}
  R_{jl}\nabla^lf=(n-1)\nabla_j R_\phi^m.
\end{equation}
Taking the covariant derivative of this equation gives
\begin{equation*}
\begin{split}
 (n-1)\nabla_i\nabla_j R_\phi^m&=R_{jl}\nabla_i\nabla^lf+\nabla_iR_{jl} \nabla_l f\\
 &=R_{ij}(\lambda-R_\phi^m)+\nabla_iR_{jl} \nabla_l f
\end{split}
\end{equation*}
where we have used the first equation in (\ref{eq2.16}) in the last equality.
Taking the trace of $i$ and $j$ and using the contracted Bianchi identity which asserts that
$\displaystyle\nabla_i R_{il}=\frac{1}{2}\nabla_l R$, we get
\begin{equation*}
  (n-1)\Delta R_\phi^m=\frac{1}{2}\langle\nabla R,\nabla f \rangle+R(\lambda-R_{\phi}^m),
\end{equation*}
which can be written as
\begin{equation}\label{eq2.27}
  (n-1)\Delta_\phi R_\phi^m=\frac{1}{2}\langle \nabla R,\nabla f\rangle+R(\lambda-R_\phi^m)-(n-1)\langle \nabla\phi,\nabla R_\phi^m\rangle.
\end{equation}
where $R$ is the scalar curvature of $g$.
Integrating \eqref{eq2.27} over $M$, we find
\begin{equation*}
  \begin{split}
  (n-1)&\int_M\Delta_\phi R_\phi^m\\
&=\frac{1}{2}\int_M \langle \nabla R,\nabla f\rangle +\int_M R(\lambda-R_\phi^m)-(n-1)\int_M\langle \nabla\phi,\nabla R_\phi^m\rangle\\
&=-\frac{1}{2}\int_M R\Delta_\phi f+\frac{1}{2}\int_{\partial M}R_\phi^m\frac{\partial f}{\partial \nu_g}+\int_M R(\lambda-R_\phi^m)-(n-1)\int_M \langle \nabla\phi,\nabla R_\phi^m\rangle\\
&=\frac{2-n-m}{2}\int_M R(\lambda-R_\phi^m)-(n-1)\int_M\langle \nabla \phi,\nabla R_\phi^m\rangle
\end{split}
\end{equation*}
where the last equality follows from \eqref{eq2.21} and \eqref{delphif}. On the other hand, by \eqref{sol2.22}, we have
\begin{equation*}
  \int_M\Delta_\phi R_\phi^m=\int_{\partial M}\frac{\partial R_\phi^m}{\partial \nu_g}=0.
\end{equation*}
Thus we obtain
\begin{equation}\label{eq2.28}
  \frac{2-n-m}{2}\int_M R(\lambda-R_\phi^m)-(n-1)\int_M\langle \nabla \phi,\nabla R_\phi^m\rangle=0.
\end{equation}

On the other hand, taking the covariant derivative of the second equation in \eqref{eq2.21}, we find
\begin{equation}\label{eq_2ndnabla}
  \nabla_j\nabla^if\nabla_i\phi+\nabla^if\nabla_j\nabla_i\phi=m\nabla_jR_\phi^m.
\end{equation}
Taking the inner product \eqref{eq_2ndnabla} with $\nabla \phi$ yields
\begin{equation*}
  \begin{split}
  m\langle \nabla R_\phi^m,\nabla \phi\rangle&=\nabla_j\nabla^if\nabla_i\phi\nabla^j\phi+\nabla^if\nabla^j\phi\nabla_j\nabla_i\phi\\
&=\nabla^2f(\nabla\phi,\nabla \phi)+\frac{1}{2}\langle\nabla f,\nabla |\nabla\phi|^2\rangle\\
&=(\lambda-R_\phi^m)|\nabla \phi|^2+\frac{1}{2}\langle\nabla f,\nabla |\nabla \phi|^2\rangle.
\end{split}
\end{equation*}
Integrating it over $M$ gives
\begin{equation}\label{delphidelr}
  \begin{split}
  \int_M\langle \nabla R_\phi^m,\nabla \phi\rangle&=\frac{1}{m}\int_M (\lambda-R_\phi^m)|\nabla \phi|^2+\frac{1}{2m}\int_M\langle \nabla f,\nabla |\nabla\phi|^2\rangle\\
&=\frac{1}{m}\int_M(\lambda-R_\phi^m)|\nabla \phi|^2-\frac{1}{2m}\int_M\Delta_\phi f|\nabla \phi|^2+\frac{1}{2m}\int_{\partial M}\frac{\partial f}{\partial \nu_g}|\nabla \phi|^2\\
&=\frac{2-n-m}{2m}\int_M(\lambda-R_\phi^m)|\nabla \phi|^2
\end{split}
\end{equation}
where the last equality follows from \eqref{eq2.21} and \eqref{delphif}.
Replacing the term $\displaystyle\int_M\langle \nabla R_\phi^m,\nabla \phi\rangle$
in \eqref{eq2.28} by
the last term in (\ref{delphidelr}), we deduce
\begin{equation}\label{eq2.23}
\int_M R(\lambda-R_\phi^m)=\frac{n-1}{m}\int_M(\lambda-R_\phi^m)|\nabla \phi|^2.
\end{equation}
By integration by parts, we have
\begin{equation}\label{eq2.24}
\begin{split}
  \int_M(\lambda-R_\phi^m)\Delta\phi=&\int_M (\lambda-R_\phi^m)(\Delta_\phi\phi+|\nabla \phi|^2)\\
=&\int_M\langle\nabla R_{\phi}^m,\nabla\phi\rangle+\int_{\partial M}(\lambda-R_\phi^m)\frac{\partial \phi}{\partial \nu_g}+\int_M(\lambda-R_\phi^m)|\nabla\phi|^2\\
=&\int_M\langle\nabla R_{\phi}^m,\nabla\phi\rangle+\int_M(\lambda-R_\phi^m)|\nabla\phi|^2
\end{split}
\end{equation}
where the last equality follows from the assumption $\frac{\partial \phi}{\partial \nu_g}=0$.
Combining (\ref{delphidelr}) and (\ref{eq2.24}), we obtain
\begin{equation}\label{eq2.25}
\int_M(\lambda-R_\phi^m)\Delta\phi
=\frac{2-n+m}{2m}\int_M(\lambda-R_\phi^m)|\nabla\phi|^2.
\end{equation}
Thus, it follows from
(\ref{eq2.23}) and (\ref{eq2.25})
that
\begin{equation}\label{int_trace2}
  \begin{split}
\int_M R_\phi^m(\lambda-R_\phi^m)
&=\int_M R(\lambda-R_\phi^m)+
2\int_M (\lambda-R_\phi^m)\Delta\phi
-\frac{m+1}{m}\int_M (\lambda-R_\phi^m)|\nabla\phi|^2\\
&=\frac{n-1}{m}\int_M(\lambda-R_\phi^m)|\nabla \phi|^2
+\frac{2-n+m}{m}\int_M(\lambda-R_\phi^m)|\nabla\phi|^2\\
&\hspace{4mm}-\frac{m+1}{m}\int_M (\lambda-R_\phi^m)|\nabla\phi|^2\\
&=0.
\end{split}
\end{equation}
Combining \eqref{int_trace} and \eqref{int_trace2}, we get
\begin{equation*}
  \int_M (\lambda-R_\phi^m)^2=0,
\end{equation*}
which implies that $R_\phi^m=\lambda$, as required.
\end{proof}

We say that the weighted Yamabe soliton with boundary \eqref{2.12def} is \textit{trivial} if the vector field $X$ is zero. It is easy to see that the trivial weighted Yamabe soliton with boundary \eqref{2.12def} has constant weighted scalar curvature in $M$. Similarly, we say that the gradient weighted Yamabe soliton with boundary \eqref{eq2.21} is \textit{trivial} if the potential function $f$ is constant. Again, it is easy to see that if the gradient weighted Yamabe soliton with boundary \eqref{eq2.21} is trivial, then it has constant weighted scalar curvature. The converse is true when $M$ is compact. To see this, integrating the first equation in (\ref{eq2.21}), we get the following:
\begin{equation}\label{cptdel}
  (n+m)\int_M(\lambda-R_\phi^m)e^{-\phi}dV_g=\int_M(\Delta_\phi f) e^{-\phi}dV_g=\int_{\partial M}\frac{\partial f}{\partial \nu_g}=0.
\end{equation}
 So  if the gradient weighted Yamabe soliton with boundary \eqref{eq2.21} has constant weighted scalar curvature $R_{\phi}^m$, then
 it follows from \eqref{cptdel} that $R_\phi^m\equiv\lambda$. Combining this with (\ref{delphif}), we have $\Delta_\phi f=0$ in $M$. Since $M$ is compact, this implies that $f$ must be constant. The following example shows that this is not the case when $M$ is not compact.

\begin{example}
 \emph{Consider $\mathbb{R}^n_+=\{x=(x_1,\cdots,x_n)\in \mathbb{R}:x_n\geq0\}$ equipped with the standard flat metric $g_0$. Then its scalar curvature and mean curvature both vanish, i.e. $R_{g_0}=0$ in $\mathbb{R}^n_+$ and $H_{g_0}=0$ on $\partial \mathbb{R}^n_+$. Take $$f(x_1,\cdots,x_n)=x_1~~\mbox{ and }~~\phi(x_1,\cdots,x_n)=a_2x_2+\cdots+a_{n-1}x_{n-1}$$
 where $a_i$ are constants and not all zero. Then one can easily see that $\nabla^2 f=0$, $\frac{\partial f}{\partial \nu_{g_0}}=0$, $\langle \nabla f,\nabla \phi\rangle =0$, $R_\phi^m=-\frac{m+1}{m}\sum_{i=2}^{n-1}a_i^2$, and $H_\phi^m=0$. Therefore, $(\mathbb{R}^n_+,g_0,f,\phi,\lambda=-\frac{m+1}{m}\sum_{i=2}^{n-1}a_i^2)$ is a nontrivial gradient weighted Yamabe soliton with boundary with constant weighted scalar curvature.}
\end{example}

\section{Blow-up Analysis on ${\mathcal{Q}}$}
In what follows, we would like to characterize Palais-Smale sequences of the functional ${\mathcal{Q}}$. In order to imitate the proof in \cite{DHR}, we consider the sequence in the following set
\begin{displaymath}
\Sigma^{+}:=\{w\in W^{1,2}(M,e^{-\phi}dV_{g})| w>0,  B(w)=1 \},
\end{displaymath}
where the functional $B(w)$ is defined by
\begin{equation}
    B(w):=\frac{\left(\int_M |w|^{p}e^{-\frac{m-1}{m}\phi}dV_g\right)^{\frac{m}{n+m-1}}}{\left(\int_{\partial M}|w|^{p}e^{-\phi}dA_{g}\right)^{\frac{n+2m-2}{n+m-1}}}, \quad p=\frac{2(m+n-1)}{n+m-2}.
\end{equation}
For general Palais-Smale sequences, the constrict $B(w)=1$ can be guaranteed by the normalization.
Equivalently, we have
\begin{equation}
    \Lambda=\inf \left\{A(w): w\in W^{1,2}(M,e^{-\phi}dV_{g}), B(w)=1\right\},
\end{equation}
where the functional $A(w)$ is defined by
\begin{equation}
    A(w):=\int_MwL_{\phi}^m(w)e^{-\phi}dV_{g}+\int_M wB_{\phi}^m(w)e^{-\phi}dA_{g}.
\end{equation}
Therefore, a Palais-Smale sequence $\{w_i\}\in \Sigma^{+}$ is characterized by the following:
\begin{center}
    the sequence ${\mathcal{Q}}(w_i)$ is bounded and \\
    $D{\mathcal{Q}}(w_i)\to 0$ strongly in $W^{1,2}(M,e^{-\phi}dV_{g})'$ as $i \to \infty$.
\end{center}
By direct calculation, it shows that $\{w_i\}$ satisfies
\begin{equation}
\begin{split}
    &D{\mathcal{Q}}(w_i). \psi=\\
    &\int_M \psi L_{\phi}^m(w_i)e^{-\phi}dV_{g}+\frac{m}{m+n-2}\frac{{\mathcal{Q}}(w_i)}{\int_M w_i^{p}e^{-\frac{m-1}{m}\phi}dV_g}\int_M \psi w_i^{p-1}e^{\frac{1}{m}\phi}e^{-\phi}dV_{g}+\\
    & \int_{\partial M} \psi B_{\phi}^m(w_i)e^{-\phi}dA_{g}-\frac{2m+n-2}{m+n-2}\frac{{\mathcal{Q}}(w_i)}{\int_{\partial M} w_i^{p}e^{-\phi}dA_{g}}\int_{\partial M} \psi w_i^{p-1}e^{-\phi}dA_{g}\\
    & \to 0
    \end{split}
\end{equation}
for all $\psi\in W^{1,2}(M,e^{-\phi}dV_{g})$; i.e.
\begin{equation}\label{E-L}
\begin{split}
    & L_{\phi}^m(w_i)+\frac{m}{m+n-2}\frac{{\mathcal{Q}}(w_i)}{\int_M w_i^{p}e^{-\frac{m-1}{m}\phi}dV_g}w_i^{p-1}e^{\frac{1}{m}\phi}\to 0 \quad \mbox{~~ in $M$~~}\\
    & B_{\phi}^m(w_i)-\frac{2m+n-2}{m+n-2}\frac{{\mathcal{Q}}(w_i)}{\int_{\partial M}w_i^{p}e^{-\phi}dA_{g}}w_i^{p-1} \to 0 \quad \mbox{~~on $\partial M$~~}
\end{split}
\end{equation}
in $W^{1,2}(M,e^{-\phi}dV_{g})'$.

\begin{theorem}\label{blowup}
Let $(M, g, e^{-\phi}dV_g,e^{-\phi}dA_g, m)$ be a smooth metric measure spaces with boundary, $m\geqslant 0$. $\{w_i\}$ is a Palais-Smale sequence in $\Sigma^{+}$ of ${\mathcal{Q}}$ such that $w_i\rightharpoonup 0$ weakly in $W^{1,2}(M,e^{-\phi}dV_{g})$ but not strongly. There exists a sequence $\{\lambda_i\}$, $\lambda_i>0$ and $\lambda_i\to +\infty$ as $i\to \infty$, a sequence $\{x_i\}$ in $M$ such that $x_i\to \tilde{x}\in \partial M$, a nontrivial nonnegative optimizer $\tilde{w} \in W^{1,2}(\R^n_+)\cap L^{\frac{2(m+n-1)}{m+n-2}}(\R^n_+)$ of the Trace Gagliardo--Nirenberg--Sobolev inequality (\ref{Trace GNS}) such that up to a subsequence, the following statement holds: if we define $w^{(2)}_i=w_i-\eta_i \tilde{w}_i$, then $\{w^{(2)}_i\}$ is also a Palais-Smale sequence of ${\mathcal{Q}}$, $w^{(2)}_i\rightharpoonup 0$ weakly in $W^{1,2}(M,e^{-\phi}dV_{g})$ and
\begin{equation}
    {\mathcal{Q}}(w^{(2)}_i)={\mathcal{Q}}(w_i)-\Lambda_{m,n}+o(1),
\end{equation}
where $\tilde{w}_i(x)=(\lambda_i)^{\frac{n-2}{2}}\tilde{w}(\lambda_i {\rm exp}^{-1}_{x_i}(x))$, $\eta_i=\eta_{\delta, x_i}$, $\delta<\frac{\textrm{inj}_g}{2}$ is a smooth cutoff function, and $||o(1)||_{W^{1,2}}\to 0$ as $i\to \infty$.
\end{theorem}
The proof of Theorem 2.1 proceeds in several steps. We follow the proof of Hebey \cite{DHR}.
\par {\bf Step 1} Palais-Smale sequences of ${\mathcal{Q}}$ in $\Sigma^{+}$ are bounded in $W^{1,2}(M,e^{-\phi}dV_{g})$.
\begin{proof}
Under the constrict $B(w_i)=1$, we know that $A(w_i)={\mathcal{Q}}(w_i)$. Since ${\mathcal{Q}}(w_i)$ is bounded and
\begin{displaymath}
    A(w)=\int_M \left(|\nabla w|^2+\frac{m+n-2}{4(m+n-1)}R^m_{\phi}w^2\right)e^{-\phi}dV_{g}+\int_{\partial M}\frac{m+n-2}{2(m+n-1)}H^m_{\phi}w^2e^{-\phi}dA_{g},
\end{displaymath}
we immediately obtain that $\{w_i\}$ is bounded in $W^{1,2}(M,e^{-\phi}dV_{g})$.
Step 1 is proved.
\end{proof}
{\bf Step 2} Subtract a Bubble $\tilde{w}_i$.
\begin{proof}
By compact embedding theorem, all lower order terms in $A(w_i)$ will vanish as $i\to \infty$. Therefore, we may assume that
\begin{displaymath}
    \lim_{i\to \infty}A(w_i)=\lim_{i\to \infty}\int_M |\nabla w_i|^2e^{-\phi}dV_{g}=\alpha\geqslant 0.
\end{displaymath}
If $\alpha=0$, we can recover the strong convergence in $W^{1,2}$. In the following, we assume that $\alpha>0$.

For $t>0$, we let
\begin{displaymath}
\mu_i(t)=\max_{x\in \overline{M}}\int_{B^+_x(t)}|\nabla w_i|^2 e^{-\phi}dV_{g},
\end{displaymath}
where $\overline{M}=M\cup \partial M$ and $B^+_x(t)=B_x(t)\cap M$.
Given $t_0>0$ small, it follows that there exists $x_0\in \overline{M}$ and $R_0>0$ such that, up to a subsequence,
\begin{equation}\label{lower bound}
\int_{B^+_{x_0}(t_0)}|\nabla w_i|^2 e^{-\phi}dV_{g}\geqslant R_0, \quad \forall i.
\end{equation}
Then, since $t\to \mu_i(t)$ is continuous, we get that for any $R\in (0, R_0)$, there exists $t_i\in (0,t_0)$ such that $\mu_i(t_i)=R$. Clearly, there also exists $x_i\in \bar{M}$ such that
\begin{displaymath}
\mu_i(t_i)=\int_{B^+_{x_i}(t_i)}|\nabla w_i|^2 e^{-\phi}dV_{g}.
\end{displaymath}
Up to a subsequence, we may assume that $x_i\to \tilde{x}$. We claim that $\tilde{x}$ is on the boundary. If it's not; i.e. $\tilde{x}\in M$, let $\hat{w}_i(x)=(\lambda_i)^{\frac{n-2}{2}}w_i(\lambda_i {\rm exp}^{-1}_{x_i}(x))$, after blow-up analysis, $\hat{w}_i(x)\to \hat{w}\geqslant 0$ strongly in $W^{1,2}_{loc}(\R^n)$. Since all lower terms in $L^m_{\phi}$ will vanish, from the first equation in (\ref{E-L}), $\hat{w}$ satisfies
\begin{equation}
    -\Delta_{\R^n} \hat{w}+c\hat{w}^{p-1}=0 \quad \mbox{~~on $\R^n$~~}
\end{equation}
for a positive constant $c$. Multiplying $\hat{w}$ and integration by parts yield that $\hat{w}=0$. It's contradicted with (\ref{lower bound}). Therefore, the atom $\tilde{x}$ is on the boundary.

We let $r_0\in (0, \frac{\textrm{inj}_g}{2})$ be such that for all $x\in \bar{M}$ and all $y,z\in \R^n$, if $|y|\leqslant r_0, |z|\leqslant r_0$, then
\begin{displaymath}
d_g({\rm exp}_x(y), {\rm exp}_x(z))\leqslant C_0|z-y|
\end{displaymath}
for some $C_0\in [1,2]$ independent of $x,y$ and $z$, where $\textrm{inj}_g$ is the injectivity radius. Given $\lambda_i\geqslant 1$ and $x\in \R^n$ such that $|x|<\textrm{inj}_g \lambda_i$, we let
\begin{displaymath}
\begin{split}
\tilde{w}_i(x) & =\lambda_i^{-\frac{n-2}{2}}{w}_i({\rm exp}_{x_i}(\lambda_i^{-1}x)),\\
\tilde{g}_i(x) & =({\rm exp}_{x_i}^{*}g)(\lambda_i^{-1}x), \\
e^{-\tilde{\phi}_i}(x) & = e^{-\phi}({\rm exp}_{x_i}(\lambda_i^{-1}x)).
\end{split}
\end{displaymath}
Then,
\begin{displaymath}
|\nabla {w}_i|_g^2({\rm exp}_{x_i}(\lambda_i^{-1}x))=\lambda^n_i|\nabla \tilde{w}_i|_{\tilde{g}_i}^2(x),
\end{displaymath}
It follows that, if $|z|+r<\textrm{inj}_g\lambda_i$, then
\begin{equation}\label{trans1}
    \int_{B^+_z(r)}|\nabla \tilde{w}_i|^2 e^{-\tilde{\phi}_i} dV_{\tilde{g}_i}=\int_{{\rm exp}_{x_i}(\lambda^{-1}B^+_z(r))}|\nabla \hat{w}_i|^2 e^{-\phi} dV_g.
\end{equation}
When $|z|+r<r_0\lambda_i$,
\begin{equation}\label{trans2}
    {\rm exp}_{x_i}(\lambda_i^{-1}B^+_z(r))\subset B^+_{{\rm exp}_{x_i}(\lambda^{-1}z)}(C_0r\lambda_i^{-1})
\end{equation}
while
\begin{equation}\label{trans3}
     {\rm exp}_{x_i}(\lambda^{-1}B^+_0(C_0r))= B^+_{x_i}(C_0r\lambda_i^{-1}).
\end{equation}
Given $r\in (0, r_0)$, we fix $t_0$ such that $C_0rt_0^{-1}\geqslant 1$. Then, for any $R\in(0, R_0),$ to be fixed later on, we let $\lambda_i\geqslant 1$ be such that $C_0r\lambda^{-1}_i=t_i$. By $(\ref{trans1})-(\ref{trans3})$, for any $z\in \R^n$ such that $|z|+r<r_0\lambda_i$,
\begin{equation}\label{R}
    \begin{split}
        \int_{B^+_z(r)}|\nabla \tilde{w}_i|^2 e^{-\tilde{\phi}_i} dV_{\tilde{g}_i} & \leqslant R \\
        \int_{B^+_0(C_0r)}|\nabla \tilde{w}_i|^2 e^{-\tilde{\phi}_i} dV_{\tilde{g}_i} &=R.
    \end{split}
\end{equation}
We let $\delta\in (0,i_g)$ and $C_1>1$ be such that for any $x\in \bar{M}$, and for any $\lambda\geqslant 1$, if $\tilde{g}_{\lambda}(y)={\rm exp}_{x}^{*}g(\lambda^{-1}y),$ $e^{-\tilde{\phi}_{\lambda}}(x)= e^{-\phi}({\rm exp}_{x}(\lambda^{-1}x))$, then
\begin{equation}\label{boundness}
\frac{1}{C_1}\int_{\R^n_+}|\nabla w|^2 dx\leqslant \int_{\R^n_+}|\nabla w|^2 e^{-\tilde{\phi}_{\lambda}} dV_{\tilde{g}_{\lambda}}\leqslant C_1\int_{\R^n_+} |\nabla w|^2 dx
\end{equation}
for all $w\in W^{1,2}(\R^n_+)$ such that $\textrm{supp}(w)\subset B^+_0(\delta \lambda)$. We let $\tilde{\eta}\in C^{\infty}_0(\R^n_+)$ be a cutoff function such that $0\leqslant \tilde{\eta}\leqslant 1$, $\tilde{\eta}=1$ in $B^+_0(1/4)$, and $\tilde{\eta}=0$ outside $B^+_0(3/4)$. We set $\tilde{\eta}_i(x)=\tilde{\eta}(\delta^{-1}\lambda_i^{-1}x)$, where $\delta$ is as above. Then
\begin{displaymath}
\int_{\R^n_+}|\nabla (\tilde{\eta}_i \tilde{w}_i)| e^{-\tilde{\phi}_{\lambda}} dV_{\tilde{g}_{\lambda}}=O(1)
\end{displaymath}
and it follows from $(\ref{boundness})$ that the sequence $\{ \tilde{\eta}_i \tilde{w}_i\}$ is bounded in $W^{1,2}(\R^n_+)$. In particular, up to a subsequence, there exists $\tilde{w}\in W^{1,2}(\R^n_+)$ such that $\tilde{\eta}_i \tilde{w}_i \rightharpoonup \tilde{w}$ weakly in $W^{1,2}(\R^n_+)$.
\end{proof}
{\bf Step 3} For $r$ and $R$ sufficiently small,
\begin{center}
    $\tilde{\eta}_i \tilde{w}_i \to \tilde{w}$ strongly in $W^{1,2}(B^+_0(C_0r))$
\end{center}
as $i\to \infty$.
\begin{proof}
We let $x_0\in \R^n_+$, and for $\rho>0$, we let $h_{\rho}$ be the standard metric on $\partial B^+_{x_0}(\rho)$. By Fatou's lemma,
\begin{displaymath}
\int^{2r}_{r}\left(\liminf_{i\to \infty}\int_{\partial B^+_{x_0}(\rho)}N_{\xi}(\tilde{\eta}_i \tilde{w}_i)dv_{h_{\rho}}\right)d\rho\leqslant \liminf_{i\to \infty}\int_{B^+_{x_0}(2r)}N_{\xi}(\tilde{\eta}_i \tilde{w}_i)dx\leqslant C
\end{displaymath}
where $N_h(w)=|\nabla w|^2+w^2$, the norm in $N_h$ is w.r.t $h$, and $\xi$ is the Euclidean metric. It follows that there exists $\rho\in [r, 2r]$ such that, up to a subsequence, and for all $i$,
\begin{displaymath}
\int_{\partial B^+_{x_0}(\rho)}N_{\xi}(\tilde{\eta}_i \tilde{w}_i)dv_{h_{\rho}} \leqslant C.
\end{displaymath}
As an easy consequence, we get that
\begin{displaymath}
||\tilde{\eta}_i \tilde{w}_i||_{W^{1,2}(\partial B^+_{x_0}(\rho))}\leqslant C
\end{displaymath}
where $C>0$ is independent of $i$. The embedding $W^{1,2}(\partial B^+_{x_0}(\rho))\subset W^{\frac{1}{2},2}(\partial B^+_{x_0}(\rho))$ is compact, and the trace operator $w\to w|_{\partial B^+}$ is continuous. It follows that, up to a subsequence, as $i\to \infty$,
\begin{displaymath}
\tilde{\eta}_i \tilde{w}_i \to \tilde{w}, \mbox{~~in~~} W^{\frac{1}{2},2}(\partial B^+_{x_0}(\rho))
\end{displaymath}
Let ${\mathcal{A}}$ be the annulus ${\mathcal{A}}=B^+_{x_0}(3r) \backslash B^+_{x_0}(\rho)$. Let also $\varphi_i\in W^{1,2}(\R^n_+)$ be such that $\varphi_i=\tilde{\eta}_i \tilde{w}_i- \tilde{w}$ in $B^+_{x_0}(\rho+\epsilon)$, and $\varphi_i=0$ in $\R^n_+ \backslash B^+_{x_0}(3r-\epsilon)$, $\epsilon>0$ small. Then
\begin{displaymath}
||\tilde{\eta}_i \tilde{w}_i- \tilde{w}||_{H^2_{1/2}(\partial B_{x_0}(\rho))}=||\varphi_i||_{H^2_{1/2}(\partial B_{x_0}(\rho))}
\end{displaymath}
while there exists $\varphi^0_i\in D^2_1({\mathcal{A}})$, the closure of $C^{\infty}_0({\mathcal{A}})$ in $W^{1,2}({\mathcal{A}})$, such that
\begin{displaymath}
||\varphi^0_i+\varphi_i||_{W^{1,2}({\mathcal{A}})}\leqslant C||\varphi_i||_{W^{\frac{1}{2},1}(\partial {\mathcal{A}})}.
\end{displaymath}
Minimization arguments give that there exits $z_i\in W^{1,2}({\mathcal{A}})$ such that
\begin{equation}
\begin{split}
    \Delta z_i=0, &\textrm{ in }{\mathcal{A}},\\
    z_i-\varphi^0_i-\varphi_i & \in D^2_1({\mathcal{A}}),
    \end{split}
\end{equation}
and
\begin{displaymath}
||z_i||_{W^{1,2}({\mathcal{A}})}\leqslant C||\varphi^0_i+\varphi_i||_{W^{1,2}({\mathcal{A}})}
\end{displaymath}
Hence, $z_i\to 0$ strongly in $W^{1,2}({\mathcal{A}})$. We let $\psi_i\in D^2_1(\R^n_+)$ be such that
\begin{displaymath}
\begin{split}
\psi_i=\left\{
  \begin{array}{ll}
    \tilde{\eta}_i \tilde{w}_i-\tilde{w}, & \hbox{$x\in \bar{B}^+_{x_0}(\rho)$;} \\
    z_i, & \hbox{$x\in \bar{B}^+_{x_0}(3r)\backslash B^+_{x_0}(\rho)$;} \\
     0, & \hbox{otherwise.}
  \end{array}
\right.
\end{split}
\end{displaymath}
We let $r$ be such that $r< \min(\textrm{inj}_g/6, \delta/24)$, and let $\tilde{\psi_i}$ be such that
\begin{displaymath}
\tilde{\psi_i}=\lambda_i^{\frac{n-2}{2}}\psi_i(\lambda_i {\rm exp}_{x_i}^{-1}(x))
\end{displaymath}
if $d_g(x_i, x)<6r$, and $\tilde{\psi_i}=0$ otherwise. Clearly, $\tilde{\eta}(\delta^{-1}{\rm exp}_{x_i}^{-1}(x))=1$ if $d_g(x_i, x)<6r$. If in addition $|x_0|<3r$, then
\begin{equation}
    D{\mathcal{Q}}(w_i). \tilde{\psi_i}=D{\mathcal{Q}}(\tilde{\eta}_i\tilde{w}_i). \psi_i.
\end{equation}
In particular, $\{\tilde{\psi}_i\}$ is bounded in $W^{1,2}(\bar{M})$. It follwos that $D{\mathcal{Q}}(w_i). \tilde{\psi_i}=o(1)$. Noting that $\psi_i \to 0$ strongly in $W^{1,2}({\mathcal{A}})$ and $\psi_i\rightharpoonup 0$ weakly in $D^2_1(\R^n_+)$, we have
\begin{equation}
    \begin{split}
    \int_{B^+_{x_0}(3r)}\langle\nabla(\tilde{\eta}_i \tilde{w}_i),\nabla \psi_i\rangle e^{-\tilde{\phi}_i} dV_{\tilde{g}_i}&=\int_{B^+_{x_0}(3r)}\langle\nabla(\psi_i+\tilde{w}),\nabla \psi_i\rangle e^{-\tilde{\phi}_i} dV_{\tilde{g}_i}\\
    &=\int_{\R^n_+}|\nabla \psi_i|^2e^{-\tilde{\phi}_i} dV_{\tilde{g}_i}+o(1).
    \end{split}
\end{equation}
Similarly, we have
\begin{equation}
    \begin{split}
        &\int_{B^+_{x_0}(3r)}(\tilde{\eta}_i \tilde{w}_i)^{p-1}\psi_i e^{-\tilde{\phi}_i} dV_{\tilde{g}_i}=\int_{\R^n_+}\psi_i^{p} e^{-\frac{m-1}{m}\tilde{\phi}_i} dV_{\tilde{g}_i}+o(1),\\
        & \int_{\partial B^+_{x_0}(3r)}(\tilde{\eta}_i \tilde{w}_i)^{p}\psi_i \tilde{v}_i e^{-\tilde{\phi}_i} dV_{\tilde{g}_i}=\int_{\partial \R^n_+}\psi_i^{p} \tilde{v}_i e^{-\tilde{\phi}_i} dV_{\tilde{g}_i}+o(1).
    \end{split}
\end{equation}
Combing them together, we obtain that
\begin{equation}\label{o(1)}
\begin{split}
    &D{\mathcal{Q}}(w_i). \tilde{\psi_i}=\\
    &\int_{\R^n_+}|\nabla \psi_i|^2e^{-\tilde{\phi}_i} dV_{\tilde{g}_i}+\frac{m}{m+n-2}\frac{{\mathcal{Q}}(w_i)}{\int_M w_i^{p}e^{-\frac{m-1}{m}\phi}dV_g}\int_{\R^n_+}\psi_i^{p} e^{-\frac{m-1}{m}\tilde{\phi}_i} dV_{\tilde{g}_i}+\\
    & -\frac{2m+n-2}{m+n-2}\frac{{\mathcal{Q}}(w_i)}{\int_{\partial M} w_i^{p}e^{-\phi}dA_{g}}\int_{\partial \R^n_+}\psi_i^{p} \tilde{v}_i e^{-\tilde{\phi}_i} dV_{\tilde{g}_i}\\
    & =o(1).
    \end{split}
\end{equation}
In the following, we would like to claim $\int_{\R^n_+}|\nabla \psi_i|^2 e^{-\tilde{\phi}_i} dV_{\tilde{g}_i}=o(1)$ from $(\ref{o(1)})$. WLOG, we assume that $\int_{\R^n_+}|\nabla \psi_i|^2 \tilde{v}_i dv_{\tilde{g}_i}\to \tilde{\alpha}$, where $\tilde{\alpha}>0$. We divide it into some cases. \\

{\bf Case 1:} Up to a subsequence, as $i\to \infty$,
\begin{equation}
\frac{\int_{\R^n_+}\psi_i^{p} e^{-\frac{m-1}{m}\tilde{\phi}_i} dV_{\tilde{g}_i}}{\int_M w_i^{p}e^{-\frac{m-1}{m}\phi}dV_g}\to +\infty \mbox{~~and~~} \frac{\int_{\partial \R^n_+}\psi_i^{p} \tilde{v}_i e^{-\tilde{\phi}_i} dV_{\tilde{g}_i}}{\int_{\partial M} w_i^{p}e^{-\phi}dA_{g}}\to +\infty.
\end{equation}
We know that if only one of them converges to infinite, we obtain the contradiction immediately. We claim that there exist positive constants $\Lambda$ and $\Gamma$, independent of $i$, such that for all $i$,
\begin{equation}\label{upper and lower}
    \gamma\leqslant \frac{\left(\int_{\R^n_+}|\nabla \psi_i|^2e^{-\tilde{\phi}_i} dV_{\tilde{g}_i}\right)\left(\int_{\R^n_+}\psi_i^{p} e^{-\frac{m-1}{m}\tilde{\phi}_i} dV_{\tilde{g}_i}\right)^{\frac{m}{m+n-1}}}{\left(\int_{\partial \R^n_+}\psi_i^{p} \tilde{v}_i e^{-\tilde{\phi}_i} dV_{\tilde{g}_i}\right)^{\frac{2m+n-2}{m+n-1}}}\leqslant \Gamma.
\end{equation}
The left inequality follows from the Trace GNS inequality and our assumption $\int_{\R^n_+}|\nabla \psi_i|^2 e^{-\tilde{\phi}_i} dV_{\tilde{g}_i}\to \tilde{\alpha}$. \\
~\\For the right one, we prove it by contradiction. If it's not true, it follows that, up to a subsequence,
\begin{displaymath}
\frac{\left(\int_{\R^n_+}\psi_i^{p} e^{-\frac{m-1}{m}\tilde{\phi}_i} dV_{\tilde{g}_i}\right)^{\frac{m}{m+n-1}}}{\left(\int_{\partial \R^n_+}\psi_i^{p} \tilde{v}_i e^{-\tilde{\phi}_i} dV_{\tilde{g}_i}\right)^{\frac{2m+n-2}{m+n-1}}}\to \infty
\end{displaymath}
Therefore, we have that for arbitrary $C>0$,
\begin{displaymath}
 C \left(\frac{\int_{\partial \R^n_+}\psi_i^{p} \tilde{v}_i e^{-\tilde{\phi}_i} dV_{\tilde{g}_i}}{\int_{\partial M} w_i^{p}e^{-\phi}dA_{g}}\right)^{\frac{2m+n-2}{m+n-1}}\leqslant \left(\frac{\int_{\R^n_+}\psi_i^{p} e^{-\frac{m-1}{m}\tilde{\phi}_i} dV_{\tilde{g}_i}}{\int_M w_i^{p}e^{-\frac{m-1}{m}\phi}dV_g}\right)^{\frac{m}{m+n-1}}.
\end{displaymath}
Choosing $C$ sufficiently large, it's contradicted with $(\ref{o(1)})$. Therefore, $(\ref{upper and lower})$ holds.\\
~\\
From $\int_{\R^n_+}|\nabla \psi_i|^2 e^{-\tilde{\phi}_i} dV_{\tilde{g}_i}\to \tilde{\alpha}$, $B(w_i)=1$ and $(\ref{upper and lower})$, we know that
\begin{equation}\label{control}
    \frac{\gamma}{\bar{\alpha}}\left(\frac{\int_{\partial \R^n_+}\psi_i^{p} \tilde{v}_i e^{-\tilde{\phi}_i} dV_{\tilde{g}_i}}{\int_{\partial M} w_i^{p}e^{-\phi}dA_{g}}\right)^{\frac{2m+n-2}{m+n-1}}\leqslant \left(\frac{\int_{\R^n_+}\psi_i^{p} e^{-\frac{m-1}{m}\tilde{\phi}_i} dV_{\tilde{g}_i}}{\int_M w_i^{p}e^{-\frac{m-1}{m}\phi}dV_g}\right)^{\frac{m}{m+n-1}} \leqslant \frac{\Gamma}{\tilde{\alpha}}\left(\frac{\int_{\partial \R^n_+}\psi_i^{p} \tilde{v}_i e^{-\tilde{\phi}_i} dV_{\tilde{g}_i}}{\int_{\partial M} w_i^{p}e^{-\phi}dA_{g}}\right)^{\frac{2m+n-2}{m+n-1}}.
\end{equation}
By our assumption and $2m+n-2>m$ when $n>2$, plugging $(\ref{control})$ into $(\ref{o(1)})$, we obtain a contradiction. {\bf Case 1} can not happen.\\

{\bf Case 2:} Up to a subsequence, as $i\to \infty$,
\begin{equation}
\frac{\int_{\R^n_+}\psi_i^{p} e^{-\frac{m-1}{m}\tilde{\phi}_i} dV_{\tilde{g}_i}}{\int_M w_i^{p}e^{-\frac{m-1}{m}\phi}dV_g}\to C_1>0 \mbox{~~and~~} \frac{\int_{\partial \R^n_+}\psi_i^{p} \tilde{v}_i e^{-\tilde{\phi}_i} dV_{\tilde{g}_i}}{\int_{\partial M} w_i^{p}e^{-\phi}dA_{g}}\to C_2>0.
\end{equation}
From {\bf Case 1:}, we know that these two limits are finite and nonnegative. Moreover, $C_1$ can't be zero. Otherwise, the ratio in $(\ref{upper and lower})$ will become as small as we want. It's contradicted with the Trace GNS inequality. Besides, $C_2$ can't be zero. If so, the left hand side of $(\ref{o(1)})$ can be bounded away from zero. It's a contradiction.\\
~\\Hence, it follows that
\begin{equation}
\frac{\left(\int_{\R^n_+}\psi_i^{p} e^{-\frac{m-1}{m}\tilde{\phi}_i} dV_{\tilde{g}_i}\right)^{\frac{m}{m+n-1}}}{\left(\int_{\partial \R^n_+}\psi_i^{p} \tilde{v}_i e^{-\tilde{\phi}_i} dV_{\tilde{g}_i}\right)^{\frac{2m+n-2}{m+n-1}}}=\frac{C_1^{\frac{m}{m+n-1}}}{C_2^{\frac{2m+n-2}{m+n-1}}}+o(1).
\end{equation}
 By strong convergence $\psi_i\to 0$ in $W^{1,2}({\mathcal{A}})$, and weak convergence $\psi_i\rightharpoonup 0$ in $D^2_1(\R^n_+)$,
\begin{displaymath}
\int_{\R^n_+}|\nabla \psi_i|^2 e^{-\tilde{\phi}_i} dV_{\tilde{g}_i}=\int_{B^+_{x_0}(\rho)}|\nabla(\tilde{\eta}_i \tilde{w}_i)|^2 e^{-\tilde{\phi}_i} dV_{\tilde{g}_i}-\int_{B^+_{x_0}(\rho)}|\nabla\tilde{w}|^2 \tilde{v}_i e^{-\tilde{\phi}_i} dV_{\tilde{g}_i}+o(1).
\end{displaymath}
so that
\begin{displaymath}
\int_{\R^n_+}|\nabla \psi_i|^2 e^{-\tilde{\phi}_i} dV_{\tilde{g}_i}\leqslant \int_{B^+_{x_0}(\rho)}|\nabla(\tilde{\eta}_i \tilde{w}_i)|^2  e^{-\tilde{\phi}_i} dV_{\tilde{g}_i}+o(1).
\end{displaymath}
Let $N$ be an integer such that $B^+_0(2)$ is covered by $N$ balls of radius $1$ and centered in $B^+_0(2)$. Then there exists $N$ points $x_i\in B^+_{x_0}(2r)$, $i=1, \cdots, N$, such that
\begin{displaymath}
B^+_{x_0}(\rho)\subset B^+_{x_0}(2r)\subset \bigcup^N_{i=1}B^+_{x_i}(r)
\end{displaymath}
and from (\ref{R}), we get that for $x_0$ and $r$ such that $|x_0|+3r<r_0$,
\begin{equation}\label{small}
    \int_{\R^n_+}|\nabla \psi_i|^2 e^{-\tilde{\phi}_i} dV_{\tilde{g}_i}\leqslant NR+o(1)
\end{equation}
Independently, we can choose $R$ sufficiently small. From $(\ref{small})$, we obtain a contradiction. Therefore, $\psi_i\to 0$ strongly in $D^2_1(\R^n_+)$. Since $r\leqslant \rho$, it follows that
\begin{center}
$\tilde{\eta}_i \tilde{w}_i\to \tilde{w}$ strongly in $W^{1,2}(B^+_{x_0}(r))$
\end{center}
and the convergence holds as soon as $NR\frac{C_1^{\frac{m}{m+n-1}}}{C_2^{\frac{2m+n-2}{m+n-1}}}< \Lambda_{m,n}$, $|x_0|<3r$, $|x_0|+3r<r_0$, $|x_0|+3r<\delta$, and $r< \min(\textrm{inj}_g/6, \delta/24)$. We fix $R, r$ sufficiently small such that $NR\frac{C_1^{\frac{m}{m+n-1}}}{C_2^{\frac{2m+n-2}{m+n-1}}}< \Lambda_{m,n}$ and $r< \min(\textrm{inj}_g/6, \delta/24)$. Then the above strong convergence holds for any $x_0$ such that $|x_0|<2r$. Since $C_0\leqslant 2$, by finite coverings, it follows that $\tilde{\eta}_i \tilde{w}_i\to \tilde{w}$ strongly in $W^{1,2}(B^+_{0}(C_0r))$. {\textbf {Step 3}} is proved.
\end{proof}
From {\bf Step 3} and $(\ref{R})$, we can write that
\begin{displaymath}
\begin{split}
    R & =\int_{B^+_{0}(C_0r)}|\nabla \tilde{w}_i|^2 e^{-\tilde{\phi}_i} dV_{\tilde{g}_i}\\
    & =\int_{B^+_{0}(C_0)}|\nabla(\tilde{\eta}_i \tilde{w}_i)|^2 e^{-\tilde{\phi}_i} dV_{\tilde{g}_i}\\
    & \leqslant C\int_{B^+_0(C_0r)}|\nabla \tilde{w}|^2 dx+o(1).
\end{split}
\end{displaymath}
It follows that $\tilde{w}\neq 0$. Let us assume that $\lambda_i\to \lambda$ as $i\to \infty$, $\lambda\geqslant 1$. If $\lambda< \infty$, then $\tilde{w}_i\rightharpoonup 0$ weakly in $W^{1,2}(B^+_0(C_0r))$ since $\hat{w}_i\rightharpoonup 0$ weakly in $H^2_1(M)$. From above, and since $\tilde{w}\neq 0$, we get that
\begin{equation}\label{scale}
\lim_{i\to \infty}\lambda_i=+\infty.
\end{equation}
Moreover, from $(\ref{small})$ and the continuity of integration, we can see that at the concentration point, the density can't vanish.

{\bf Step 4} For any $\lambda>0$,
\begin{displaymath}
    \tilde{w}_i\to \tilde{w} \mbox{~~strongly in~~} W^{1,2}(B^+_0(\lambda))
\end{displaymath}
as $i\to \infty$. Moreover, $\tilde{w}$ satisfies the equation
\begin{equation}\label{E-L eq}
    \begin{split}
        & -\Delta \tilde{w}+c_1\tilde{w}^{p-1}e^{\frac{\phi}{m}}(\tilde{x})=0 \quad \mbox{~~in $\R^n_+$~~}\\
        & \frac{\partial }{\partial \nu}\tilde{w}-c_2\tilde{w}^{p-1}=0 \quad \mbox{~~ on $\partial \R^n_+$~~},
    \end{split}
\end{equation}
where
\begin{displaymath}
\begin{split}
    & c_1=\frac{m}{m+n-2}\frac{{\mathcal{Q}[\R^n_+, dx^2]}(\tilde{w})}{\int_{\R^n_+} \tilde{w}^{p}e^{-\frac{m-1}{m}\phi}(\tilde{x})dx}, \\
    & c_2=\frac{2m+n-2}{m+n-2}\frac{{\mathcal{Q}[\R^n_+, dx^2]}(\tilde{w})}{\int_{\partial \R^n_+} \tilde{w}^{p}e^{-\phi}(\tilde{x})dx}.
\end{split}
\end{displaymath}
\begin{proof}
We let $\lambda\geqslant 1$ be given. By $(\ref{scale})$, for $i$ large, $\lambda_i\geqslant \lambda$, and $(\ref{R})$ holds for $z$ such that $|z|+r<r_0\lambda$. Then, as is easily checked from the proof of {\bf Step 3}, local strong convergence holds if $|x_0|<3r(2\lambda-1)$, $|x_0|+3r<r_0\lambda$, and $|x_0|+3r<\delta\lambda$. In particular, local strong convergence holds if $|x_0|<2r\lambda.$ Hence, $\tilde{\eta}_i \tilde{w}_i\to \tilde{w}$ strongly in $W^{1,2}(B^+_0(2r\lambda))$. Noting that for $x$ in a compact subset of $\R^n_+$, $\tilde{\eta}_i(x)=1$ for $i$ large enough, and that $\lambda\geqslant 1$ is arbitrary,we easily get that $\tilde{w}_i\to \tilde{w}$ strongly in $W^{1,2}(B^+_0(\lambda))$.

Next, we prove the second argument. Let $\varphi\in C^{\infty}_0(\R^n_+)$ be the test function, and $\lambda_0>0$ be such that $\textrm{supp}(\varphi)\subset B^+_0(\lambda_0)$. Let also $\hat{\varphi}_i$ be given by
\begin{displaymath}
\hat{\varphi}_i(x)=\lambda_i^{\frac{n-2}{2}}\varphi(\lambda_i x).
\end{displaymath}
Then $\textrm{supp}(\hat{\varphi}_i)\subset B^+_0(\lambda_i^{-1}\lambda_0)$. For $i$ large, we let $\varphi_i$ be the smooth function on $\overline{M}$ given by $\hat{\varphi}_i=\varphi_i\circ {\rm exp}_{x_i}$. Notice that $\tilde{g}_i\to \xi$ and $e^{-\tilde{\phi}_i}(x)\to e^{-\phi}(\tilde{x})$ in $C^1(B^+_0(\lambda))$ for arbitrary $\lambda>0$. For $i$ large, direct computation shows that
\begin{equation}
    \int_M (\nabla w_i \nabla \varphi_i)e^{-\phi} dV_g=\int_{\R^n_+}\langle\nabla (\tilde{\eta}_i \tilde{w}_i),\nabla \varphi\rangle e^{\tilde{\phi}_i} dV_{\tilde{g}_i}=\int_{\R^n_+}\langle\nabla \tilde{w}, \nabla \varphi\rangle e^{-\phi}(\tilde{x})dx+o(1).
\end{equation}
By the strong convergence in $W^{1,2}(\R^n)$ and the compact embedding theorem, all lower order terms will vanish. Hence, taking the limit in (\ref{o(1)}), the second argument is proved.
\end{proof}

For $x\in \overline{M}$ and $\hat{\delta}\in (0, \delta/8)$, we let $v_i$ be given by
\begin{equation}
    \tilde{w}_i(x)=\eta_i(x)\lambda_i^{\frac{n-2}{2}}\tilde{w}(\lambda_i {\rm exp}_{x_i}^{-1}(x))
\end{equation}
where $\eta_i=\eta_{\hat{\delta}, x_i}$. We let $w^{(2)}_i=w_i-\eta_i \tilde{w}_i$ and claim the following holds.

{\bf Step 5} The following relations hold. On the one hand,
\begin{equation}\label{weak cong}
w^{(2)}_i\rightharpoonup 0
\end{equation}\label{energy}
weakly in $W^{1,2}(M)$, as $i\to \infty$. On the other hand,
\begin{equation}\label{small2}
 D{\mathcal{Q}}(\tilde{w}_i)\to 0 \mbox{~~and~~} D{\mathcal{Q}}(w^{(2)}_i)\to 0
\end{equation}
strongly as $i\to \infty$. Finally,
\begin{equation}
     {{\mathcal{Q}}}(w^{(2)}_i)= {{\mathcal{Q}}}(w_i)-\Lambda_{m,n}+o(1).
\end{equation}
\begin{proof}
We start with the proof of $(\ref{weak cong})$. It suffices to prove that $\tilde{w}_i\rightharpoonup 0$ weakly in $W^{1,2}(\overline{M})$. This weak convergence follows from  the blow-up property of $\tilde{w}_i$, since we have $\lambda_i\to +\infty$, as $i\to \infty$. This proves $(\ref{weak cong})$. Besides, by embedding theorems, we obtain strong convergence in $L^s$, $s<p$. It proves $(\ref{small2})$. Finally, for $(\ref{energy})$, we can lift them to Euclidean half space $\R^n_+$, and similar argument in {\bf Step 3} implies
\begin{equation}
    {{\mathcal{Q}}}(w^{(2)}_i)= {{\mathcal{Q}}}(w_i)-{{\mathcal{Q}}[\R^n_+,dx^2]}(\tilde{w})+o(1).
\end{equation}
Notice that $\tilde{w}$ satisfies the equation (\ref{E-L eq}). Therefore, $\tilde{w}$ is an optimizer of the Trace GNS inequality (\ref{Trace GNS}) with ${{\mathcal{Q}}[\R^n_+,dx^2]}(\tilde{w})=\Lambda_{m,n}$.
{\bf Step 5} is proved.
\end{proof}

\section{Aubin-type almost trace GNS inequalities}
In this section, we adapt the argument in \cite[Theorem 1.1]{Yan} to obtain almost trace GNS inequalities. Combining this with Theorem \ref{blowup} yields the proof of Posso's conjecture.

\begin{lemma}\label{AB}
    Let $(M^n, \partial M,g)$ be a compact Riemannian manifold with boundary. Then there exists positive constants $A,B$ depending on $(M^n, \partial M,g)$ such that for all $w\in C^{\infty}(M)$,
    \begin{equation}
    \begin{split}
        \left(\int_{\partial M}|w|^{\frac{2(n+m-1)}{n+m-2}}dA_{g}\right)^{\frac{n+2m-2}{n+m-1}} \leqslant &\left(A\int_M |\nabla w|^2 dV_{g}+B\int_M w^2e^{-\phi}dV_{g}\right)\\
        \times &\left(\int_M |w|^{\frac{2(n+m-1)}{n+m-2}}dV_g\right)^{\frac{m}{n+m-1}}.
    \end{split}
    \end{equation}
\end{lemma}
\begin{proof}
    Fix $\epsilon>0$. Around each point $P\in \partial M$ we can choose a neighborhood $U\subset \overline{M}$ such that, in normal coordinates on $U$, the eigenvalues of metric $g$ are between $(1+\epsilon)^{-1}$ and $(1+\epsilon)$, and furthermore $(1+\epsilon)^{-1}dx\leqslant dA_g\leqslant (1+\epsilon)dx$. Choose a finite subcover $\{U_i\}_{i=1}^k$ of $\partial M$ and a subordinate partition of unity, which we may write as $\{\alpha^2_i\}_{i=1}^k$, where $\alpha_i\in C^{\infty}(M)$ and $\sum \alpha_i^2=1$. Then we have
    \begin{equation}\label{subadd}
        \begin{split}
            \left(\int_{\partial M}|w|^{\frac{2(n+m-1)}{n+m-2}}dA_{g}\right)^{\frac{n+m-2}{2(n+m-1)}}&=\left(\int_{\partial M}\left(\sum \alpha^2_i w^2\right)^{\frac{n+m-1}{n+m-2}}dA_{g}\right)^{\frac{n+m-2}{2(n+m-1)}}\\
            & \leqslant \sum \left(\int_{U_i}|\alpha_i w|^{\frac{2(n+m-1)}{n+m-2}}dA_{g}\right)^{\frac{n+m-2}{n+m-1}}\\
            & \leqslant C(\epsilon) \sum \left(\int_{U_i}|\alpha_i w|^{\frac{2(n+m-1)}{n+m-2}}dx\right)^{\frac{n+m-2}{n+m-1}},
        \end{split}
    \end{equation}
    for some constant $C(\epsilon)$ depending on $\epsilon$.
    Applying the Trace GNS inequality (\ref{Trace GNS}) on each $U_i$ yields
    \begin{equation}
        \left(\int_{U_i}|\alpha_i w|^{\frac{2(n+m-1)}{n+m-2}}dx\right)^{\frac{n+2m-2}{n+m-1}}\leqslant \frac{1}{\Lambda_{m,n}}\left(\int_{U_i}|\nabla (\alpha_iw)|_{0}^2dx\right)\left(\int_{U_i} |\alpha_i w|^{\frac{2(m+n-1)}{m+n-2}}dx\right)^{\frac{m}{m+n-1}}.
    \end{equation}
    By estimates on the deviation of $g$ from the Euclidean metric and Cauchy--Schwartz inequality, there exists constants $C_1, C_2$ depending on $\epsilon$ and $\{\alpha_i\}$ such that
    \begin{equation}
        \int_{U_i}|\nabla (\alpha_iw)|_{0}^2dx\leqslant C_1 \int_{U_i}\alpha_i^2|\nabla w|_{g}^2dV_g+C_2\int_{U_i}w^2dV_g.
    \end{equation}
    Therefore, we obtain that
    \begin{equation}\label{upp}
        \begin{split}
            &\left(\int_{U_i}|\alpha_i w|^{\frac{2(n+m-1)}{n+m-2}}dx\right)^{\frac{n+2m-2}{n+m-1}}\\
            &\leqslant \frac{1}{\Lambda_{m,n}} \left(C_1 \int_{U_i}\alpha_i^2|\nabla w|_{g}^2dV_g+C_2\int_{U_i}w^2dV_g\right)\left(\int_{U_i} |\alpha_i w|^{\frac{2(m+n-1)}{m+n-2}}dx\right)^{\frac{m}{m+n-1}}\\
            &\leqslant \frac{1}{\Lambda_{m,n}} \left(C_1 \int_{M}\alpha_i^2|\nabla w|_{g}^2dV_g+C_2\int_{M}w^2dV_g\right)\left(\int_{M} |w|^{\frac{2(m+n-1)}{m+n-2}}dx\right)^{\frac{m}{m+n-1}}.
        \end{split}
    \end{equation}
    Combing the uniform upper bound in (\ref{upp}) with the subadditivity (\ref{subadd}) yields the desired result.
\end{proof}

\begin{theorem}\label{Aubin}
    Let $(M, g, e^{-\phi}dV_g,e^{-\phi}dA_g, m)$ be a smooth metric measure spaces with boundary. For any $\epsilon>0$, there exists a constant $C(\epsilon)$ such that for all $w\in  W^{1,2}(M,e^{-\phi}dV_{g})$, we have
    \begin{equation}
        \begin{split}
        &\left(\int_{\partial M}|w|^{\frac{2(n+m-1)}{n+m-2}}e^{-\phi}dA_{g}\right)^{\frac{n+2m-2}{n+m-1}} \leqslant \left(\int_M |w|^{\frac{2(n+m-1)}{n+m-2}}e^{-\frac{m-1}{m}\phi}dV_g\right)^{\frac{m}{n+m-1}}\\
        & \times \left(\left(\frac{1}{\Lambda_{m,n}}+\epsilon\right)\int_M |\nabla w|^2e^{-\phi}dV_{g}
        +C(\epsilon)\left(\int_M w^2e^{-\phi}dV_{g}+\int_{\partial M}w^2e^{-\phi}dA_{g}\right)\right).
        \end{split}
    \end{equation}
\end{theorem}
\begin{proof}
    Let $\alpha=\frac{1}{\Lambda_{m,n}}+\epsilon$. If it's not true, for any $i\in \N$, we can find a function $w_i$ in $W^{1,2}(M,e^{-\phi}dV_{g})$ with $B(w_i)=1$ such that
    \begin{equation}\label{contr}
        1>\alpha \int_M |\nabla w_i|^2e^{-\phi}dV_{g}+ i \left(\int_M w_i^2e^{-\phi}dV_{g}+\int_{\partial M}w_i^2e^{-\phi}dA_{g}\right).
    \end{equation}
    It implies that $\{w_i\}$ is bounded in $W^{1,2}(M,e^{-\phi}dV_{g})$. Therefore, we obtain that $w_i\rightharpoonup 0$ in $W^{1,2}(M,e^{-\phi}dV_{g})$ and $w_i\to 0$ in $L^{2}(M,e^{-\phi}dV_{g})$. According to Lemma \ref{AB} and $B(w_i)=1$, we conclude that $\int_M |\nabla w_i|^2e^{-\phi}dV_{g}$ is bounded away from zero uniformly. Repeating Theorem \ref{blowup}, we obtain that
    \begin{equation}
        w_i-\sum_{j=1}^{N}\eta_i^j\tilde{w}^j_i\to 0
    \end{equation}
    strongly in $W^{1,2}(M,e^{-\phi}dV_{g})$ and for each $j$, $\tilde{w}^j_i\to \tilde{w}^j$, an optimizer of (\ref{Trace GNS}), $N\in \N, N\geqslant 1$. Besides, we have ${{\mathcal{Q}}}(w_i)=N\Lambda_{m,n}+o(1)$. Therefore, from
    \begin{equation}
        N\Lambda_{m,n}<\frac{1}{\alpha}=\frac{\Lambda_{m,n}}{\Lambda_{m,n}\epsilon+1}<\Lambda_{m,n}.
    \end{equation}
    It's a contradiction.
\end{proof}
\begin{proof}[Proof of Theorem \ref{main}]
    According to \cite[Theorem A]{Posso18}, it suffices to consider the case $0\leqslant \Lambda<\Lambda_{m,n}$. Let $\{w_i\}$ be a minimizing sequence. Without loss of generality, we assume that $B(w_i)=1$ and $w_i\rightharpoonup w$ weakly in $W^{1,2}(M,e^{-\phi}dV_{g})$. If $w$ is zero, by Theorem \ref{Aubin}, we obtain that
    \begin{displaymath}
        1\leqslant \left(\frac{1}{\Lambda_{m,n}}+\epsilon\right) \int_M |\nabla w_i|^2e^{-\phi}dV_{g}<1.
    \end{displaymath}
    It's a contradiction. Therefore, $w$ is nonzero. By compact embedding theorem and Fatou's lemma, we get that
    \begin{equation}
        A(w)\leqslant \Lambda, \quad B(w)=1.
    \end{equation}
    It implies that $w$ is a minimizer. By regularity argument mentioned in \cite{Posso18}, we know that $w$ is a smooth positive minimizer.
\end{proof}

\section{Acknowledgement}

The first author was supported
by Basic Science Research Program through the National Research Foundation of Korea (NRF) funded by the Ministry of Education, Science and Technology (Grant No. 2020R1A6A1A03047877 and 2019R1F1A1041021), and by Korea Institute for Advanced Study (KIAS) grant
funded by the Korea government (MSIP). The second author was supported by a KIAS Individual Grant (SP070701) via the Center for Mathematical Challenges at Korea Institute for Advanced Study.

\bibliographystyle{amsplain}

\begin{thebibliography}{30}

\bibitem{ADN}
S. Agmon, A. Douglis, and L. Nirenberg,
Estimates near the boundary for solutions of elliptic partial differential equations satisfying general boundary conditions. I.
\textit{Comm. Pure Appl. Math.} \textbf{12} (1959), 623–727.



\bibitem{Aubin0} T. Aubin, \'{E}quations diff\'{e}rentielles non lin\'{e}aires et probl\`{e}me de Yamabe concernant
la courbure scalaire. \textit{J. Math. Pures Appl. (9)} \textbf{55} (1976), 269-296.





\bibitem{Case15}
J. S. Case,
A Yamabe-type problem on smooth metric measure spaces. \textit{J. Differential Geom.} \textbf{101} (2015), no. 3, 467–505.


\bibitem{Case12}
J. S. Case,
Conformal invariants measuring the best constants for Gagliardo-Nirenberg-Sobolev inequalities.
\textit{Calc. Var. Partial Differential Equations} \textbf{48} (2013), no. 3-4, 507–526.


\bibitem{Case13}
J. S. Case,
Sharp metric obstructions for quasi-Einstein metrics.
\textit{J. Geom. Phys.} \textbf{64} (2013), 12–30.



\bibitem{Case19}
J. S. Case,
The weighted $\sigma_k$-curvature of a smooth metric measure space. \textit{Pacific J. Math.} \textbf{299} (2019), no. 2, 339–399.


\bibitem{Souza20}
M. de Souza,
On the existence of extremals for the weighted Yamabe problem on compact manifolds.
\textit{Differential Geom. Appl.} \textbf{68} (2020), 101585, 16 pp.




\bibitem{Escobar2}
J. F. Escobar,
Conformal deformation of a Riemannian metric to a scalar flat metric with constant mean curvature on the boundary. \textit{Ann. of Math. (2)} \textbf{136} (1992), 1-50.

\bibitem{Escobar1}
J. F. Escobar,  The Yamabe problem on manifolds with boundary. \textit{J. Differential Geom.} \textbf{35} (1992), 21-84.

\bibitem{Escobar3}
J. F. Escobar,
Uniqueness and non-uniqueness of metrics with prescribed scalar and mean curvature on compact manifolds with boundary. \textit{J. Funct. Anal.} \textbf{202} (2003), no. 2, 424–442.

\bibitem{Lee} J.M. Lee, Introduction to Smooth Manifolds. Springer, New York (2003)

\bibitem{Han}
Y. Han and H. Lin,  Vanishing theorems for $f$-harmonic forms on smooth metric measure spaces. \textit{Nonlinear Anal.} \textbf{162} (2017), 113–127.



\bibitem{Nirenberg}
L. Nirenberg,   Topics in nonlinear functional analysis. Chapter 6 by E. Zehnder. Notes by R. A. Artino. Revised reprint of the 1974 original. \textit{Courant Lecture Notes in Mathematics}, \textbf{6}. New York University, Courant Institute of Mathematical Sciences, New York; American Mathematical Society, Providence, RI, 2001. xii+145 pp.


\bibitem{Posso21}
J. M. Posso, A generalization of Aubin's result for a Yamabe-type problem on smooth metric measure spaces. \textit{Bull. Sci. Math.} \textbf{172} (2021), Paper No. 103052, 33 pp.



\bibitem{Posso18}
J. M. Posso, A generalization of Escobar-Riemann mapping type problem on smooth metric measure spaces.
\textit{Comm. Anal. Geom.} (2019), accepted. https://arxiv.org/abs/1805.03694


\bibitem{S}
M. Schechter,   On $L^p$ estimates and regularity. II. \textit{Math. Scand.} \textbf{13} (1963), 47–69.

\bibitem{Schoen} R. Schoen, Conformal deformation of a Riemannian metric to constant scalar curvature.
\textit{J. Differential Geom.} \textbf{20} (1984), 479-495.


\bibitem{Trudinger} N. S. Trudinger, Remarks concerning the conformal deformation of Riemannian structures
on compact manifolds. \textit{Ann. Scuola Norm. Sup. Pisa (3)} \textbf{22} (1968), 265-274.



\bibitem{Yamabe} H. Yamabe, On a deformation of Riemannian structures on compact manifolds. \textit{Osaka
Math. J.} \textbf{12} (1960), 21-37.

\bibitem{BCFG} F. Borell, D. Cordero--Erausquin, Y. Fujita, I. Gentil, A. Guillin, New Sharp Gagliardo–Nirenberg–Sobolev Inequalities and an Improved Borell–Brascamp–Lieb Inequality. \textit{IMRN} \textbf{10} (2020), 3042–3083.

\bibitem{DHR} O. Druet, E. Hebey, F. Robert, Blow-up Theory for Elliptic PDEs in Riemannian Geometry.

\bibitem{Yan} Z. Yan, Improved Higher-order Sobolev inequalities. \textit{arXiv}: 2203.13873.
\end{thebibliography}

\end{document}